\documentclass[11pt,a4paper]{amsart}

\pdfoutput=1

\usepackage{amssymb}
\usepackage{amsmath}
\usepackage{amsthm}
\usepackage{amsbsy}
\usepackage{enumerate}
\usepackage{xcolor}
\usepackage{esint}
\usepackage{comment}
\usepackage{graphicx}

\usepackage[T1]{fontenc}
\usepackage[utf8]{inputenc}
\usepackage[english]{babel}

\theoremstyle{plain}
\newtheorem{theorem}{Theorem}[section]
\newtheorem{proposition}[theorem]{Proposition}
\newtheorem{lemma}[theorem]{Lemma}

\theoremstyle{definition}

\newtheorem{remark}[theorem]{Remark}
\theoremstyle{remark}
\newtheorem*{discussion}{Discussion}

\numberwithin{equation}{section}

\newcommand*{\E}{\mathbb{E}}
\newcommand*{\R}{\mathbb{R}}

\newcommand*{\Var}{\operatorname{Var}}

\newcommand*{\EE}{\mathbb E}
\newcommand*{\PP}{\mathbb P}

\newcommand*{\RR}{\mathbb R}
\newcommand*{\NN}{\mathbb{N}}

\newcommand*{\bbN}{\mathbb N}

\newcommand*{\bbR}{\mathbb R}

\newcommand*{\cA}{\mathcal A}

\newcommand*{\cG}{\mathcal G}

\newcommand*{\cI}{\mathcal I}

\newcommand*{\cL}{\mathcal L}

\DeclareMathOperator{\const}{const}
\DeclareMathOperator{\sgn}{sgn}

\renewcommand*{\doteq}{:=}

\begin{document}
\title{Truncated control variates for weak approximation schemes}\thanks{The work of Denis Belomestny is supported by the Russian Science Foundation project 14-50-00150.}

\author{Denis Belomestny}
\address{University of Duisburg-Essen, Essen, Germany and IITP RAS, Moscow, Russia}
\email{denis.belomestny@uni-due.de}

\author{Stefan H\"afner}
\address{PricewaterhouseCoopers GmbH, Frankfurt, Germany}
\email{stefan.haefner@de.pwc.com}

\author{Mikhail Urusov}
\address{University of Duisburg-Essen, Essen, Germany}
\email{mikhail.urusov@uni-due.de}


\begin{abstract}
In this paper we present an enhancement of the regression-based variance reduction approaches recently proposed in Belomestny et al.\ \cite{belomestny2016variance} and~\cite{belomestny2016stratified}.
This enhancement is based on a truncation of the control variate and allows for a significant reduction of the computing time, while the complexity stays of the same order.
The performances of the proposed truncated algorithms are illustrated by a numerical example.


\end{abstract}

\maketitle

\section*{Introduction}
\label{sec:intro}
Let \(T>0\) be a fixed time horizon.
Consider a $d$-dimensional diffusion process
$(X_t)_{t\in[0,T]}$
defined by the It\^o stochastic differential equation
\begin{align}\label{x_sde}
dX_t
=\mu(X_t)\,dt
+\sigma(X_t)\,dW_{t},
\quad X_{0}=x_0\in\mathbb{R}^d,
\end{align}
for Lipschitz continuous functions
\(\mu\colon\mathbb{R}^d\to\mathbb{R}^d\)
and
\(\sigma\colon\mathbb{R}^d
\to\mathbb{R}^{d\times m}\),
where \((W_t)_{t\in[0,T]}\)
is a standard \(m\)-dimensional Brownian motion.
Our aim is to compute the expectation 
\begin{align}\label{eq:2512a0}
u(t,x)\doteq\E[ f(X_{T}^{t,x})],
\end{align}
for some \(f:\) \(\mathbb{R}^d\to \mathbb{R},\) where $X^{t,x}$ denotes the solution to~\eqref{x_sde} started
at time $t$ in point~$x$. The standard Monte Carlo (SMC)
estimate for \(u(0,x)\) at a fixed point \(x\in \bbR^d\) has the form
\begin{align}\label{eq:2512a1}
V_{N_0}\doteq\frac{1}{N_0}\sum_{i=1}^{N_0} f\Bigl(\overline{X}_T^{(i)}\Bigr)
\end{align} 
for some \(N_0\in \mathbb{N}_0,\) where \(\overline{X}_T\) is an approximation for \(X^{0,x}_T\) constructed via a time discretisation of~\eqref{x_sde} (we refer to~\cite{KP} for a nice overview of various discretisation schemes).
In the computation of $u(0,x)=\EE [f(X^{0,x}_T)]$
by the SMC approach
there are two types of error inherent:
the (deterministic) discretisation error
$\EE [f(X^{0,x}_T)]-\EE [f(\overline X_{T})]$
and the Monte Carlo (statistical) error,
which results from the substitution of
$\EE [f(\overline{X}_{T})]$
with the sample average~$V_{N_0}$.
The aim of variance reduction methods is to reduce the latter statistical error.  For example, in the so-called control variate variance reduction approach
one looks for a random variable
\(\xi\) with \(\EE \xi=0\),
which can be simulated,
such that
the variance of the difference
\(f(\overline{X}_{T})-\xi\) is minimised, that~is,
\begin{align*}
\Var[f(\overline{X}_{T})-\xi]\to\min\text{ under } \EE\xi=0.
\end{align*}
Then one uses the sample average
\begin{align}\label{eq:2512a2}
V_{N_0}^{CV}\doteq\frac{1}{N_0}\sum_{i=1}^{N_0}
\left[f\Bigl(\overline{X}_T^{(i)}\Bigr)-\xi^{(i)}\right]
\end{align}
instead of~\eqref{eq:2512a1} to approximate $\EE [f(\overline{X}_{T})]$.
The use of control variates for computing expectations of functionals of diffusion processes via Monte Carlo  was initiated by Newton~\cite{newton1994variance} and further developed in Milstein and Tretyakov~\cite{milstein2009practical}.
In Belomestny et al~\cite{belomestny2016variance} a novel regression-based approach for the construction of control variates, which reduces the variance of  the approximated functional \(f(\overline{X}_{T})\) was proposed. 
As shown in~\cite{belomestny2016variance},
the ``Monte Carlo approach with
the Regression-based Control Variate''
(abbreviated below as ``RCV approach'') is able to achieve a higher order convergence
of  the resulting variance to zero,
which in turn leads to a significant complexity reduction
as compared to the SMC algorithm. Other prominent examples of algorithms with this property are the multilevel Monte Carlo (MLMC) algorithm of~\cite{giles2008multilevel} and quadrature-based algorithms of~\cite{muller2015complexity} and~\cite{muller1324deterministic}.
The RCV approach becomes especially simple
in the case of the so-called weak approximation schemes,
i.e.,\ the schemes, where simple random variables
are used in place of Brownian increments,
and which became quite popular in recent years. However, due to the fact that a lot of computations are required for implementing the RCV approach, its numerical efficiency is not convincing in higher-dimensional examples.
The same applies also to the SRCV algorithm of~\cite{belomestny2016stratified}.
In this paper we further enhance the performances
of the RCV and SRCV algorithms by truncating the control variates, leading to a reduction from $(2^m-1)$ to $m$ terms at each time point in case of the weak Euler scheme and a reduction from $(3^m2^\frac{m(m-1)}{2}-1)$ to $m(m+1)=O(m^2)$ terms at each time point in case of the second order weak scheme. It turns out that, while the computing time is reduced significantly, we still have a sufficient variance reduction effect such that the complexity is of the same order as for the original RCV and SRCV approaches. 

The paper is organised as follows.
In Section~\ref{sec:discr} we present a smoothness theorem for a general class of discretisation schemes.
Section~\ref{sec:weak} recalls the construction of control variates for weak schemes of the first and the second order.
The main truncation results are derived in Section~\ref{sec:trunc}.
In Section~\ref{sec:regr} we describe a generic regression algorithm. Section~\ref{sec:compl} deals with a complexity analysis for the algorithm that is based on the truncated control variate.
Section~\ref{sec:num} is devoted to a simulation study.
Finally, all proofs are collected in Section~\ref{sec:proofs}.

\section{Smoothness theorem for discretisation schemes}
\label{sec:discr}
In this section we present a technical result for discretisation schemes, which will be very important in the sequel.
To begin with, let
$J\in\bbN$ denote the time discretisation parameter,
we set $\Delta\doteq T/J$
and consider discretisation schemes
defined on the grid
$\{j\Delta:j=0,\ldots,J\}$. 

Let us consider a scheme,
where $d$-dimensional approximations
$X_{\Delta,j\Delta}$,
$j=0,\ldots,J$, satisfy $X_{\Delta,0}=x_0$ and
\begin{align}
\label{smooth:eq:defX}
    X_{\Delta, j\Delta} 
    = 
    \Phi_{\Delta}\left(
       X_{\Delta, (j-1)\Delta},\xi_j\right),\quad j=1,\ldots,J,
\end{align}
for some Borel measurable functions
$\Phi_{\Delta}\colon\bbR^{d+ \tilde m}\to\bbR^d$, where $\tilde m\ge m$, and for $\tilde m$-dimensional
i.i.d. random vectors $\xi_j=(\xi_j^1,\ldots,\xi_j^{\tilde m})^\top$ with independent coordinates satisfying $\EE\left[\xi_j^i\right]=0$ and $\Var\left[\xi_j^i\right]=1$ for all $i=1,\ldots,\tilde m$, $j=1,\ldots,J$. Moreover, let $\cG_0$ be the trivial $\sigma$-field
and $\cG_j=\sigma(\xi_1,\ldots,\xi_j)$, $j=1,\ldots,J$. In the chapters below we will focus on different kinds of discretisation schemes, resulting in different convergence behaviour.

We now define the random function $G_{l,j}(x)$ for $J\ge l\ge j\ge 0$, $x\in\RR^d$, as follows
\begin{align}
\label{smooth:def_G}
G_{l,j}(x)&\equiv\Phi_{\Delta,l}\circ\Phi_{\Delta,l-1}\circ\ldots\circ
\Phi_{\Delta,{j+1}}(x),\quad l>j,\\
\notag
G_{l,j}(x)&\equiv x,\quad l=j,
\end{align}
where $\Phi_{\Delta,l}(x)\doteq\Phi_{\Delta}\left(x,\xi_l\right)$ for $l=1,\ldots,J$. By $\Phi_{\Delta,l}^k$, $k\in\left\{1,\ldots,d\right\}$, we denote the $k$-th component of the function $\Phi_{\Delta,l}$. 
Note that 
it holds 
\begin{align}
\label{smooth:rel_uG}
q_j(x):=\EE\left[f(X_{\Delta,T})\left|X_{\Delta,j\Delta}\right.=x\right]=\EE\left[f(G_{J,j}(x))\right].
\end{align}
Let us define the operator $D^\alpha$ as follows
\begin{align}
\label{smooth:operat_D}
D^\alpha g(x)
&\doteq
\frac{\partial^{\left|\alpha\right|}g(x)}{\partial x_1^{\alpha_1}\cdots\partial x_d^{\alpha_{d}}},
\end{align}
where $g$ is a real-valued function, $\alpha\in\NN_0^{d}$ and $\left|\alpha\right|=\alpha_1+\ldots+\alpha_d$
($\NN_0:=\NN\cup\{0\}$).

In the next theorem we present some smoothness conditions on $q_j$, which will be used several times in the chapters below.
\begin{theorem}
\label{smooth:func_assump_general}
Let $K\in\left\{1,2,3\right\}$. Suppose that $f$ is $K$ times continuously differentiable with bounded partial derivatives up to order $K$,
$\Phi_\Delta(\cdot,\xi)$ is $K$ times continuously
differentiable (for any fixed~$\xi$),
and that, for any $n\in\NN$, $l\ge j$, $k\in\left\{1,\ldots,d\right\}$, $\alpha\in\NN_0^d$ with $1\le\left|\alpha\right|\le K$, it holds
\begin{align}
\label{smooth:Phi_assump_1}
\left|\EE\left[\left.\left(D^\alpha\Phi{}_{\Delta,l+1}^{k}(G_{l,j}(x))\right)^{n}\right|\mathcal{G}_{l}\right]\right|\leq\begin{cases}
(1+A_{n}\Delta), & \left|\alpha\right|=\alpha_k=1\\
B_{n}\Delta, & (\left|\alpha\right|>1)\vee (\alpha_k\ne 1)
\end{cases}
\end{align}
with probability one for some constants $A_{n}>0$, $B_{n}>0$.
Moreover, suppose that for any $n_1,n_2\in\NN$, $\alpha,\beta\in\NN_0^d$,
with $|\alpha|=1$,
$1\le\left|\beta\right|\le K$,
$\alpha\ne\beta$,
it holds
\begin{align}
\label{smooth:Phi_assump_2 new}
\left|\EE\left[\left.\left(D^{\alpha}\Phi{}_{\Delta,l+1}^{k}(G_{l,j}(x))\right)^{n_1}\left(D^{\beta}\Phi{}_{\Delta,l+1}^{k}(G_{l,j}(x))\right)^{n_2}\right|\mathcal{G}_{l}\right]\right|\leq E_{n_1,n_2}\Delta
\end{align}
for some constants $E_{n_1,n_2}>0$. Then we obtain for all $j\in\left\{0,\ldots,J\right\}$ that $q_j$ is $K$ times continuously differentiable with bounded partial derivatives up to order~$K$. 
\end{theorem}

\section{Representations for weak approximation schemes}
\label{sec:weak}

Below we focus on weak schemes of first and second order.

\subsection{Weak Euler scheme}
\label{sec:first}
In this subsection we treat weak schemes of order~$1$.
Let us consider a scheme,
where $d$-dimensional approximations
$X_{\Delta,j\Delta}$,
$j=0,\ldots,J$, satisfy $X_{\Delta,0}=x_0$ and
\begin{align}
\label{eq:scheme_structure_md}
X_{\Delta,j\Delta}=
\Phi_{\Delta}(X_{\Delta,(j-1)\Delta},\xi_j),
\quad j=1,\ldots,J,
\end{align}
for some functions
$\Phi_{\Delta}\colon\bbR^{d+m}\to\bbR^d$,
with $\xi_j=(\xi_j^1,\ldots,\xi_j^m)$, $j=1,\ldots,J$,
being $m$-dimensional
iid random vectors with iid coordinates
such that
\begin{align*}
\PP\left(\xi_j^k=\pm1\right)=\frac12, \quad k=1,\ldots,m.
\end{align*}
That is, relating to the framework in Section~\ref{sec:discr}, we have $\tilde m=m$ and use the discrete increments $\xi_j^i$, $i=1,\ldots,m$.
A particular case is the weak Euler scheme (also called the
\emph{simplified weak Euler scheme}
in \cite[Section~14.1]{KP})
of order~1, which is given by
\begin{align}
\label{eq:PhiK=1}
\Phi_{\Delta}(x,y)
=x+\mu(x)\,\Delta+\sigma(x)\,y\,\sqrt{\Delta}.
\end{align}
Let us recall the functions (cf.~\eqref{smooth:rel_uG})
\begin{align*}
q_j(x)=\EE[f(X_{\Delta, T})|X_{\Delta,j\Delta}=x].
\end{align*}
The proposition below summarises important representations for the weak Euler scheme, which were derived in~\cite{belomestny2016variance}.
\begin{proposition}
\label{th:weak_md01}
The following representation holds
\begin{align}
\label{eq:repr02}
f(X_{\Delta,T})=\EE f(X_{\Delta,T})
+\sum_{j=1}^J
\sum_{r=1}^m
\sum_{1\le s_1<\ldots<s_r\le m}
a_{j,r,s}(X_{\Delta,(j-1)\Delta})
\prod_{i=1}^r \xi_j^{s_i},
\end{align}
where we use the notation
$s=(s_1,\ldots,s_r)$.
The coefficients
$a_{j,r,s}\colon\bbR^d\to\bbR$
can be computed by the formula
\begin{align}
\label{eq:coef05}
a_{j,r,s}(x)
=\EE\left[\left.
f(X_{\Delta,T}) \prod_{i=1}^r \xi_j^{s_i} 
\,\right|\,
X_{\Delta,(j-1)\Delta}=x
\right]
\end{align}
for all $j$, $r$, and $s$ as in~\eqref{eq:repr02}. Moreover, we have the following recursion formulas
\begin{align}
\notag
q_{j-1}(x)=&\EE\bigl[q_j(X_{\Delta, j\Delta})|X_{\Delta,(j-1)\Delta}=x\bigr]
=\frac{1}{2^{m}}\sum_{y=(y^{1},\ldots,y^{m})\in\left\{ -1,1\right\} ^{m}}q_j(\Phi_\Delta(x,y)),\\
\label{eq:coef05a}
a_{j,r,s}(x)=&\frac{1}{2^m}\sum_{y=(y^{1},\ldots,y^{m})\in\left\{ -1,1\right\} ^{m}}\,  \left [\prod_{i=1}^r y^{s_i}\right]
q_j(\Phi_\Delta(x,y)),
\end{align}
for all $j\in\left\{1,\ldots,J\right\}$, $r\in\left\{1,\ldots,m\right\}$, $1\le s_1<\ldots<s_r\le m$, where $q_J\equiv f$.
\end{proposition}

The next proposition (cf. Proposition~3.2 in~\cite{belomestny2016variance}) shows the properties of the weak Euler scheme combined with the control variate 
\begin{align}
\label{eq:2909a2}
M^{(1)}_{\Delta,T}\doteq
\sum_{j=1}^J
\sum_{r=1}^m
\sum_{1\le s_1<\ldots<s_r\le m}
a_{j,r,s}(X_{\Delta,(j-1)\Delta})
\prod_{i=1}^r \xi_j^{s_i}, 
\end{align}
where the coefficients
$a_{j,r,s}(x)$ are given by~\eqref{eq:coef05}.

\begin{proposition}
\label{prop:Euler:CV}
Assume that $\mu$ and $\sigma$ in~\eqref{x_sde} are Lipschitz continuous with components $\mu^i,\,\sigma^{i,r}\colon \R^d\to\R$,
$i=1,\ldots,d$, $r=1,\ldots,m$,
being $4$ times continuously differentiable
with their partial derivatives of order up to $4$
having polynomial growth.
Let $f\colon\R^d\to\R$ be $4$ times continuously differentiable with 
partial derivatives of order up to $4$
having polynomial growth.
Provided that \eqref{eq:PhiK=1} holds
and that, for sufficiently large $p\in\mathbb N$,
the expectations $\EE |X_{\Delta,j\Delta}|^{2p}$
are uniformly bounded in $J$ and $j=0,\ldots,J$,
we have for this
``simplified weak Euler scheme''
\begin{align*}
\left|\E\left[f(X_T) - f(X_{\Delta,T})\right]\right|\le c\Delta,
\end{align*}
where the constant $c$ does not depend on $\Delta$. Moreover, it holds
$\Var\left[f(X_{\Delta,T}) - M^{(1)}_{\Delta,T}\right]=0.$
\end{proposition}

\begin{discussion}
In order to use the control variate
$M^{(1)}_{\Delta,T}$
in practice, we need to estimate the unknown coefficients
$a_{j,r,s}$.
Thus, practically implementable control variates
$\widetilde{M}^{(1)}_{\Delta,T}$
have the form~\eqref{eq:2909a2}
with some estimated functions
$\tilde{a}_{j,r,s}\colon\bbR^d\to\bbR$.
Notice that they remain valid control variates,
i.e.\ we still have $\EE\bigl[\widetilde{M}^{(1)}_{\Delta,T}\bigr]=0$,
which is due to the martingale transform
structure\footnote{\label{ft:19102016a1}This phrase means
that the discrete-time process
$\tilde M=(\tilde M_l)_{l=0,\ldots,J}$, where $\tilde M_0=0$ and
$\tilde M_l$ is defined like the right-hand side of~\eqref{eq:2909a2}
but with
$\sum_{j=1}^J$ being replaced by $\sum_{j=1}^l$
and $a_{j,r,s}$ by $\tilde a_{j,r,s}$
is a martingale, which is a straightforward calculation.}
in~\eqref{eq:2909a2}.
\end{discussion}

\subsection{Second order weak scheme}
\label{sec:second}

Now we treat weak schemes of order~$2$.
We consider a scheme, where
$d$-dimensional approximations
$X_{\Delta,j\Delta}$, $j=0,\ldots,J$, satisfy
$X_{\Delta,0}=x_0$ and
\begin{align}
\label{eq:2002a5}
X_{\Delta,j\Delta}=
\Phi_{\Delta}(X_{\Delta,(j-1)\Delta},\xi_j,V_j),
\quad j=1,\ldots,J,
\end{align}
for some functions
$\Phi_{\Delta}\colon\bbR^{d+m+m\times m}\to\bbR^d$.
Here,
\begin{itemize}
\item[(S1)]
$\xi_j=(\xi_j^i)_{i=1}^m$
are $m$-dimensional random vectors,
\item[(S2)]
$V_j=(V_j^{il})_{i,l=1}^m$
are random $m\times m$-matrices,
\item[(S3)]
the pairs $(\xi_j,V_j)$, $j=1,\ldots,J$, are i.i.d.,
\item[(S4)]
for each $j$, the random elements $\xi_j$ and $V_j$
are independent,
\item[(S5)]
for each $j$, the random variables
$\xi_j^i$, $i=1,\ldots,m$, are i.i.d.\ with
\begin{align*}
\PP\left(\xi_j^i=\pm\sqrt{3}\right)=\frac16,
\quad
\PP\left(\xi_j^i=0\right)=\frac23,
\end{align*}
\item[(S6)]
for each $j$, the random variables
$V_j^{il}$, $1\le i<l\le m$, are i.i.d.\ with
\begin{align*}
\PP\left(V_j^{il}=\pm1\right)=\frac12,
\end{align*}
\item[(S7)]
$V_j^{li}=-V_j^{il}$, $1\le i<l\le m$, $j=1,\ldots,J$,
\item[(S8)]
$V_j^{ii}=-1$, $i=1,\ldots,m$, $j=1,\ldots,J$.
\end{itemize}
Hence, the matrices $V_j$ can be generated by means of $\frac{m(m-1)}{2}$ i.i.d. random variables. That is, relating to the framework in Section~\ref{sec:discr}, we have $\tilde m$-dimensional random vectors $\tilde\xi_j:=((\xi_j^i)_{i=1,\ldots,m},(V_j^{il})_{1\le i<l\le m})$ with $\tilde m=m+\frac{m(m-1)}{2}=\frac{m(m+1)}{2}$.

\begin{remark}
In order to obtain a second order weak scheme
in the multidimensional case,
we need to incorporate additional
random elements $V_j$
into the structure of the scheme.
This is the reason
why we now consider~\eqref{eq:2002a5}
instead of~\eqref{eq:scheme_structure_md}.
For instance, to get the second order weak scheme (also called the
\emph{simplified order~2 weak Taylor scheme})
of \cite[Section~14.2]{KP}
in the multidimensional case,
we need to define the functions
$\Phi_{\Delta}(x,y,z)$,
$x\in\bbR^d$, $y\in\bbR^m$, $z\in\bbR^{m\times m}$,
as explained below.
First we define the function
$\Sigma\colon\bbR^d\to\bbR^{d\times d}$
by the formula
\begin{align*}
\Sigma(x)=\sigma(x)\sigma(x)^\top
\end{align*}
and recall that the coordinates
of vectors and matrices are denoted
by superscripts,
e.g.\ $\Sigma(x)=(\Sigma^{kl}(x))_{k,l=1}^d$,
$\Phi_{\Delta}(x,y,z)
=(\Phi_{\Delta}^k(x,y,z))_{k=1}^d$.
Let us introduce the operators
$\cL^r$, $r=0,\ldots,m$,
that act on sufficiently smooth functions
$g\colon\bbR^d\to\bbR$ as follows:
\begin{align*}
\cL^0 g(x)&\doteq\sum_{k=1}^d
\mu^k(x) \frac{\partial g}{\partial x^k}(x)
+\frac12 \sum_{k,l=1}^d
\Sigma^{kl}(x) \frac{\partial^2 g}{\partial x^l\partial x^k}(x),
\\
\cL^r g(x)&\doteq\sum_{k=1}^d \sigma^{kr}(x)
\frac{\partial g}{\partial x^k}(x),\quad
r=1,\ldots,m. 
\end{align*}
The $r$-th coordinate $\Phi_{\Delta}^r$,
$r=1,\ldots,d$, in the simplified order~2
weak Taylor scheme of \cite[Section~14.2]{KP}
is now given by the formula
\begin{align}
\Phi_{\Delta}^r(x,y,z)&=
x^r+\sum_{k=1}^m \sigma^{rk}(x)\,y^k\,\sqrt{\Delta}
\label{eq:2002a6}\\
&\hspace{1em}+\left[
\mu^r(x)+\frac12\sum_{k,l=1}^m
\cL^k\sigma^{rl}(x) (y^k y^l+z^{kl})
\right]\Delta
\notag\\
&\hspace{1em}+\frac12\sum_{k=1}^m
\left[
\cL^0\sigma^{rk}(x)+\cL^k \mu^r(x)
\right]
y^k\,\Delta^{3/2}
+\frac12\cL^0\mu^r(x)\,\Delta^2,
\notag
\end{align}
provided the coefficients $\mu$ and $\sigma$
of~\eqref{x_sde}
are sufficiently smooth.
We will need to work
explicitly with~\eqref{eq:2002a6}
at some point,
but all results in this subsection
assume structure~\eqref{eq:2002a5} only.
\end{remark}

Let us define the index sets
\begin{align*}
\cI_1=\{1,\ldots,m\},\quad
\cI_2=\left\{(k,l)\in\cI_1^2:k<l\right\}
\end{align*}
and the system
\begin{align*}
\cA=\left\{(U_1,U_2)\in\mathcal P(\cI_1)\times\mathcal
P(\cI_2):U_1\cup U_2\ne\emptyset\right\},
\end{align*}
where $\mathcal P(\cI)$ denotes
the set of all subsets of a set~$\cI$.
For any $U_1\subseteq\cI_1$
and $o\in\{1,2\}^{U_1}$,
we write $o$ as
$o=(o_r)_{r\in U_1}$.
Below we use the convention
that a product over the empty set
is always one.

For $k\in\bbN_0$,
$H_k\colon\bbR\to\bbR$
stands for the (normalized) $k$-th Hermite polynomial,~i.e.
\begin{align*}
  H_k(x) 
  \doteq 
  \frac{(-1)^k}
    {\sqrt{k!}}
  e^{\frac{x^2}{2}}
  \frac{d^k}{dx^k}e^{-\frac{x^2}{2}},
  \quad x\in\bbR.
\end{align*}
We remark that, in particular,
$H_0\equiv1$, $H_1(x)=x$ and $H_2(x)=\frac1{\sqrt2}(x^2-1)$.

As in Subsection~\ref{sec:first}, we summarise important representations from~\cite{belomestny2016variance} below. 

\begin{proposition}
\label{th:weak_md03}
It holds
\begin{align}
\label{eq:2002a1}
f(X_{\Delta,T})=\EE f(X_{\Delta,T})
+\sum_{j=1}^J
\sum_{(U_1,U_2)\in\cA}
\sum_{o\in\{1,2\}^{U_1}}
a_{j,o,U_1,U_2}(X_{\Delta,(j-1)\Delta})
\prod_{r\in U_1} H_{o_r}(\xi_j^r)
\prod_{(k,l)\in U_2} V_j^{kl},
\end{align}
where the coefficients
$a_{j,o,U_1,U_2}\colon\bbR^d\to\bbR$
can be computed by the formula
\begin{align}
\label{eq:2002a2}
a_{j,o,U_1,U_2}(x)
=\EE\left[\left.
f(X_{\Delta,T})
\prod_{r\in U_1} H_{o_r}(\xi_j^r)
\prod_{(k,l)\in U_2} V_j^{kl}
\right| X_{\Delta,(j-1)\Delta}=x
\right].
\end{align}
Moreover, we have for each $j\in\{1,\ldots,J\}$,
\begin{align*}
q_{j-1}(x)=&\EE[q_{j}(X_{\Delta,j\Delta})|X_{\Delta,(j-1)\Delta}=x]
\\
=&\frac{1}{2^{\frac{m(m-1)}{2}}}\,\frac{1}{6^{m}}\sum_{(y^{1},\ldots,y^{m})\in\{-\sqrt{3},0,\sqrt{3}\}^{m}}\sum_{(z^{uv})_{1\le u<v\le m}\in\{-1,1\}^{\frac{m(m-1)}{2}}}4^{\sum_{i=1}^{m}I(y^{i}=0)}q_j(\Phi_\Delta(x,y,z)),
\end{align*}
and, for all $(U_1,U_2)\in\cA$, $o\in\{1,2\}^{U_1}$,
it holds
\begin{align}
\label{eq:0403a11}
a_{j,o,U_1,U_2}(x)=\frac{1}{2^{\frac{m(m-1)}{2}}}\,\frac{1}{6^{m}}\sum_{(y^{1},\ldots,y^{m})\in\{-\sqrt{3},0,\sqrt{3}\}^{m}}\sum_{(z^{uv})_{1\le u<v\le m}\in\{-1,1\}^{\frac{m(m-1)}{2}}}\\
\notag 4^{\sum_{i=1}^{m}I(y^{i}=0)}
\prod_{r\in U_1}H_{o_r}(y^r)
\prod_{(k,l)\in U_2}z^{kl}\,
q_j(\Phi_\Delta(x,y,z)),
\end{align}
where $y=(y^1,\ldots,y^m)$, $z=(z^{uv})$ is the $m\times m$-matrix
with $z^{vu}=-z^{uv}$, $u<v$, $z^{uu}=-1$ and $q_J\equiv f$.
\end{proposition}

Using Theorem~\ref{th:weak_md03},
we obtain the following result (see Proposition~3.6 in~\cite{belomestny2016variance}),
which provides a bound for the discretisation error
and a perfect control variate for the discretised quantity.

\begin{proposition}
\label{prop:second:CV}
Assume, that $\mu$ and $\sigma$ in~\eqref{x_sde} are Lipschitz continuous with components
$\mu^i,\,\sigma^{i,r}\colon \R^d\to\R$,
$i=1,\ldots,d$, $r=1,\ldots,m$,
being $6$ times continuously differentiable
with their partial derivatives of order up to $6$
having polynomial growth.
Let $f\colon\R^d\to\R$ be $6$ times continuously differentiable
with partial derivatives of order up to $6$
having polynomial growth.
Provided that~\eqref{eq:2002a6} holds
and that, for sufficiently large $p\in\mathbb N$,
the expectations $\EE |X_{\Delta,j\Delta}|^{2p}$
are uniformly bounded in $J$ and $j=0,\ldots,J$,
we have for this
``simplified second order weak Taylor scheme''
\begin{align*}
\left|\E\left[f(X_T) - f(X_{\Delta,T})\right]\right|\le c\Delta^2,
\end{align*}
where the constant $c$ does not depend on $\Delta$. Moreover, we have
$\Var\left[f(X_{\Delta,T}) - M^{(2)}_{\Delta,T}\right]=0$
for the control variate
\begin{align}
\label{eq:28042016a1}
M^{(2)}_{\Delta,T}\doteq\sum_{j=1}^J
\sum_{(U_1,U_2)\in\cA}
\sum_{o\in\{1,2\}^{U_1}}
a_{j,o,U_1,U_2}(X_{\Delta,(j-1)\Delta})
\prod_{r\in U_1} H_{o_r}(\xi_j^r)
\prod_{(k,l)\in U_2} V_j^{kl},
\end{align}
where the coefficients
$a_{j,o,U_1,U_2}(x)$ are defined in~\eqref{eq:2002a2}.
\end{proposition}

%
%
%

\section{Truncated control variates for weak approximation schemes}
\label{sec:trunc}
Below we recall the assumptions from~\cite{belomestny2016variance},
suggest sufficient conditions for them in terms of the functions
$f,\mu,\sigma$,
and then suggest some stronger conditions
that will justify the use of truncated control variates.

\subsection{Weak Euler scheme}
\label{sec:trunc_first}
Note that we considered only the second order weak scheme in terms of the regression and complexity analyses in~\cite{belomestny2016variance}. However, analogous assumptions for the weak Euler scheme are as follows (cf. Proposition~\ref{th:weak_md01}): fix some $j\in\left\{1,\ldots, J\right\}$, $r\in\left\{1,\ldots,m\right\}$, $s=(s_1,\ldots,s_r)$ with $1\le s_1<\ldots<s_r\le m$, set $\zeta_{j,r,s}\doteq f(X_{\Delta,T})\prod_{i=1}^r \xi_j^{s_i}$
and remark that
$a_{j,r,s}(x)=\EE[\zeta_{j,r,s}|X_{\Delta, (j-1)\Delta}=x
]$. We assume that, for some positive
constants $\Sigma,A$, it holds:
\begin{itemize}
\item[(A1)]
$\sup_{x\in\R^d}\Var[\zeta_{j,r,s}|X_{\Delta,(j-1)\Delta}=x]
\le\Sigma<\infty$,
\item[(A2)]
$\sup_{x\in\R^d} |a_{j,r,s}(x)|\le A\sqrt{\Delta}<\infty$.
\end{itemize}
In the following theorem we suggest sufficient conditions
for the above assumptions.

\begin{theorem}
\label{th:first_old}
(i)
Let $f$ be bounded. Then (A1) holds.

(ii) Let all the functions $\sigma^{ki}$, $k\in\left\{1,\ldots,d\right\}$, $i\in\left\{1,\ldots,m\right\}$, be bounded and all the functions $f,\mu^k,\sigma^{ki}$ be continuously differentiable with bounded partial derivatives. Then (A2) holds.
\end{theorem}

Next we suggest some stronger conditions
that give us somewhat more than~(A2).

\begin{theorem}
\label{th:first_new}
Let all the functions $\sigma^{ki}$, $k\in\left\{1,\ldots,d\right\}$, $i\in\left\{1,\ldots,m\right\}$, be bounded and all the functions $f,\mu^k,\sigma^{ki}$ be twice continuously differentiable with bounded partial derivatives up to order 2. Then it holds
\begin{itemize}
\item[(A3)] $\sup_{x\in\R^d} |a_{j,r,s}(x)|\lesssim\Delta$,
whenever $r>1$.
\end{itemize}
\end{theorem}

\begin{remark}
As a generalisation of Theorem~\ref{th:first_new}, it is natural to expect that it holds, under additional smoothness conditions on $f,\mu,\sigma$,
$$
\sup_{x\in\RR^d}\left|a_{j,r,s}(x)\right|\lesssim \Delta^{r/2}
$$
for all $j\in\left\{1,\ldots,J\right\}$, $r\in\left\{1,\ldots,m\right\}$ and $1\le s_1< \ldots < s_r\le m$.
\end{remark}

Let us define the ``truncated control variate'' 
\begin{align}
\label{cv:first_trunc}
M_{\Delta,T}^{(1),trunc}\doteq \sum_{j=1}^J\sum_{i=1}^ma_{j,1,e_i}(X_{\Delta,(j-1)\Delta})\xi_j^i,
\end{align}
where $e_i\in\RR^m$ denotes the $i$-th unit vector in $\RR^m$ and $a_{j,1,e_i}$ is given by (cf.~\eqref{eq:coef05})
\begin{align*}
a_{j,1,e_i}(x)
=\EE\left[\left.
f(X_{\Delta,T})\xi_j^i 
\,\right|\,
X_{\Delta,(j-1)\Delta}=x
\right].
\end{align*}
Note that the superscript ``trunc'' comes from ``truncated''. That is, we consider in $M_{\Delta,T}^{(1),trunc}$ only the terms of the control variate $M_{\Delta,T}^{(1)}$ for which $r=1$ (cf.~\eqref{eq:2909a2}). 

Next we study the truncation error that arises from replacing $M_{\Delta,T}^{(1)}$ by $M_{\Delta,T}^{(1),trunc}$.
\begin{theorem}
\label{th:first_var}
Suppose that all the functions $\sigma^{ki}$, $k\in\left\{1,\ldots,d\right\}$, $i\in\left\{1,\ldots,m\right\}$ are bounded and all the functions $f,\mu^k,\sigma^{ki}$ are twice continuously differentiable with bounded partial derivatives up to order 2. Then it holds (cf.\ Proposition~\ref{prop:Euler:CV})
\begin{align}
\label{first:cv_var}
\Var\left[f(X_{\Delta,T})-M_{\Delta,T}^{(1),trunc}\right]\lesssim\Delta.
\end{align}
\end{theorem}
Notice that under Assumption~(A2) alone
the variance in~\eqref{first:cv_var} would have been $O(1)$.

\subsection{Second order weak scheme}
\label{sec:trunc:second}
First we recall some of the required assumptions in~\cite{belomestny2016variance}:
let us fix some $j\in\left\{1,\ldots, J\right\}$, $(U_1,U_2)\in\mathcal{A}$, $o\in\left\{1,2\right\}^{U_1}$, set
\begin{align*}
\zeta_{j,o,U_1,U_2}\doteq f(X_{\Delta,T})
\prod_{r\in U_1} H_{o_r}(\xi_j^r)
\prod_{(k,l)\in U_2} V_j^{kl}
\end{align*}
and remark that
$a_{j,o,U_1,U_2}(x)=\EE[\zeta_{j,o,U_1,U_2}|X_{\Delta, (j-1)\Delta}=x
]$. We assume that, for some positive
constants $\Sigma,A$, it holds:
\begin{itemize}
\item[(B1)]
$\sup_{x\in\R^d}\Var[\zeta_{j,o,U_1,U_2}|X_{\Delta,(j-1)\Delta}=x]
\le\Sigma<\infty$,
\item[(B2)]
$\sup_{x\in\R^d} |a_{j,o,U_1,U_2}(x)|\le A\sqrt{\Delta}<\infty$.
\end{itemize}
Below we verify the above assumptions.

\begin{theorem}
\label{th:second_old}
(i)
Let $f$ be bounded. Then (B1) holds.

(ii) Let all the functions $\mu^k$ and $\sigma^{ki}$, $k\in\left\{1,\ldots,d\right\}$, $i\in\left\{1,\ldots,m\right\}$, be bounded,
the function $f$ be continuously differentiable
with bounded partial derivatives,
and all the functions $\mu^k,\sigma^{ki}$ be three times
continuously differentiable with bounded partial derivatives up to order~3. Then (B2) holds.
\end{theorem}

Let us define the index sets
$$
\mathcal{K}_1:=\left\{r\in U_1:o_r=1\right\},\quad
\mathcal{K}_2:=\left\{r\in U_1:o_r=2\right\}.
$$
In the following theorem we provide some stronger
conditions that give us more than~(B2).

\begin{theorem}
\label{th:second_new}
(i)
Let all the functions $\mu^k$ and $\sigma^{ki}$, $k\in\left\{1,\ldots,d\right\}$, $i\in\left\{1,\ldots,m\right\}$, be bounded,
the function $f$ be twice continuously differentiable
with bounded partial derivatives up to order~2,
and all the functions $\mu^k,\sigma^{ki}$ be four times continuously differentiable with bounded partial derivatives up to order~4. Then it holds
\begin{itemize}
\item[(B3)] $\sup_{x\in\R^d} |a_{j,o,U_1,U_2}(x)|\lesssim\Delta$, whenever $\left|U_2\right|+\left|\mathcal{K}_2\right|+\frac{\left|\mathcal{K}_1\right|}{2}\ge 1$.
\end{itemize}

(ii)
Let all the functions $\mu^k$ and $\sigma^{ki}$, $k\in\left\{1,\ldots,d\right\}$, $i\in\left\{1,\ldots,m\right\}$, be bounded,
the function $f$ be three times continuously differentiable
with bounded partial derivatives up to order~3,
and all the functions $\mu^k,\sigma^{ki}$ be five times continuously differentiable with bounded partial derivatives up to order~5. Then it holds
\begin{itemize}
\item[(B4)] $\sup_{x\in\R^d} |a_{j,o,U_1,U_2}(x)|\lesssim\Delta^{3/2}$, whenever $\left|U_2\right|+\left|\mathcal{K}_2\right|+\frac{\left|\mathcal{K}_1\right|}{2}>1$.
\end{itemize}
\end{theorem}

\begin{remark}
(i) As a generalisation of Theorem~\ref{th:second_new}, it is natural to expect that it holds, under additional smoothness conditions on $f,\mu,\sigma$,
$$
\sup_{x\in\RR^d}\left|a_{j,o,U_1,U_2}(x)\right|\lesssim \Delta^{\left|U_2\right|+\left|\mathcal{K}_2\right|+\frac{\left|\mathcal{K}_1\right|}{2}}
$$
for all $j\in\left\{1,\ldots,J\right\}$, $(U_1,U_2)\in\cA$ and $o\in\{1,2\}^{U_1}$.

(ii) Define 
\begin{align}
\label{delta_u}
\Delta_{U_1,U_2}:=\left\{\begin{array}{ll} \Delta^{\left|U_2\right|+\left|\mathcal{K}_2\right|+\frac{\left|\mathcal{K}_1\right|}{2}} & \text{if }\left|U_2\right|+\left|\mathcal{K}_2\right|+\frac{\left|\mathcal{K}_1\right|}{2}\le 1, \\
\Delta^{3/2} & \text{otherwise}.\end{array}\right.
\end{align}
An equivalent reformulation of assumptions~(B2)--(B4) is as follows: there exists some positive constant $\tilde A$ such that it holds
\begin{align}
\label{tilde_A}
\sup_{x\in\RR^d}\left|a_{j,o,U_1,U_2}(x)\right|\le \tilde A\,\Delta_{U_1,U_2}
\end{align}
for all $j,o,U_1,U_2$.
\end{remark}
Similar to Section~\ref{sec:trunc_first}, let us define a truncated control variate through
\begin{align}
\label{cv:second_trunc}
M_{\Delta,T}^{(2),trunc}\doteq \sum_{j=1}^J
\sum_{\substack{(U_1,U_2)\in\cA\\ \left|U_2\right|+\left|\mathcal{K}_2\right|+\frac{1}{2}\left|\mathcal{K}_1\right|\le 1}}
\sum_{o\in\{1,2\}^{U_1}}
a_{j,o,U_1,U_2}(X_{\Delta,(j-1)\Delta})
\prod_{r\in U_1} H_{o_r}(\xi_j^r)
\prod_{(k,l)\in U_2} V_j^{kl}.
\end{align}

Next we derive the truncation error that arises from replacing $M_{\Delta,T}^{(2)}$ by $M_{\Delta,T}^{(2),trunc}$.

\begin{theorem}
\label{th:second_var}
Suppose that all the functions $\mu^k$ and $\sigma^{ki}$, $k\in\left\{1,\ldots,d\right\}$, $i\in\left\{1,\ldots,m\right\}$ are bounded,
the function $f$ is three times continuously differentiable
with bounded partial derivatives up to order~3,
and all the functions $\mu^k,\sigma^{ki}$ are five times continuously differentiable with bounded partial derivatives up to order~5. Then it holds (cf.\ Proposition~\ref{prop:second:CV})
\begin{align}
\label{second:cv_var}
\Var\left[f(X_{\Delta,T})-M_{\Delta,T}^{(2),trunc}\right]\lesssim\Delta^2.
\end{align}
\end{theorem}

\section{Generic regression algorithm}
\label{sec:regr}
In the previous sections
we have given several representations
for control variates.
Now we discuss how to compute the coefficients
in these representations via regression.
For the sake of clarity,
we focus on second order schemes
and control variate~\eqref{cv:second_trunc}
with coefficients given by~\eqref{eq:2002a2}.

\subsection{Monte Carlo regression}
Fix a $Q$-dimensional vector of real-valued functions \ensuremath{\psi=(\psi^{1},\ldots,\psi^{Q})}
on \ensuremath{\mathbb{R}^{d}}. Simulate
a big number\footnote{In the complexity analysis below
we show how large $N$ is required to be
in order to provide an estimate within some
given tolerance.}
$N$ of independent  ``training paths'' of the
discretised diffusion $X_{\Delta,j\Delta},$ $j=0,\ldots, J$.
In what follows these $N$ training paths
are denoted by $D_N^{tr}$:
$$
D_N^{tr}\doteq
\left\{
(X_{\Delta,j\Delta}^{tr,(i)})_{j=0,\ldots,J}:
i=1,\ldots,N
\right\}.
$$
Let
$\boldsymbol{\alpha}_{j,o,U_1,U_2}=(
\alpha_{j,o,U_1,U_2}^{1},\ldots,\alpha_{j,o,U_1,U_2}^{Q})$,
where $j\in\left\{1,\ldots, J\right\}$, $(U_1,U_2)\in\mathcal{A}$, $\left|U_2\right|+\left|\mathcal{K}_2\right|+\frac{1}{2}\left|\mathcal{K}_1\right|\le 1$, $o\in\left\{1,2\right\}^{U_1}$,
 be a solution of the following least squares optimisation problem:
\begin{align*}
\operatorname{argmin}_{\boldsymbol{\alpha}\in\mathbb{R}^{n}}
\sum_{i=1}^{N}\left[\zeta^{tr,(i)}_{j,o,U_1,U_2}-\alpha^{1}\psi^{1}(X_{\Delta,(j-1)\Delta}^{tr,(i)})-\ldots-\alpha^{Q}\psi^{Q}(X_{\Delta,(j-1)\Delta}^{tr,(i)})\right]^{2}
\end{align*}
with 
\begin{align*}
\zeta_{j,o,U_1,U_2}^{tr,(i)}\doteq f(X_{\Delta,T}^{tr,(i)})
\prod_{r\in U_1} H_{o_r}\left((\xi_j^{tr,(i)})^r\right)
\prod_{(k,l)\in U_2} (V_j^{tr,(i)})^{kl}.
\end{align*}
Define an estimate for  the coefficient function $a_{j,o,U_1,U_2}$ via
\begin{align*}
\hat a_{j,o,U_1,U_2}(x)\doteq
\hat a_{j,o,U_1,U_2}(x,D_N^{tr})\doteq
\alpha_{j,o,U_1,U_2}^{1}\psi^{1}(x)+\ldots+\alpha_{j,o,U_1,U_2}^{Q}\psi^{Q}(x),\quad x\in\mathbb{R}^{d}.
\end{align*}
The intermediate expression
$\hat a_{j,o,U_1,U_2}(x,D_N^{tr})$
in the above formula
emphasises that the estimates
$\hat a_{j,o,U_1,U_2}$
of the functions $a_{j,o,U_1,U_2}$
are random in that they depend on
the simulated training paths.
The cost of computing
$\boldsymbol{\alpha}_{j,o,U_1,U_2}$ is of order $O(NQ^{2})$,
 since each \ensuremath{\boldsymbol{\alpha}_{j,o,U_1,U_2}}
 is of the form \ensuremath{\boldsymbol{\alpha}_{j,o,U_1,U_2}=B^{-1}b}
 with 
\begin{align}
\label{b_matr_reg}
B_{k,l}\doteq\frac{1}{N}\sum_{i=1}^{N}\psi^{k}\bigl(X_{\Delta,(j-1)\Delta}^{tr,(i)}\bigr)\psi^{l}\bigl(X_{\Delta,(j-1)\Delta}^{tr,(i)}\bigr)
\end{align}
and 
\begin{align*}
b_{k}\doteq\frac{1}{N}\sum_{i=1}^{N}\psi^{k}\bigl(X_{\Delta,(j-1)\Delta}^{tr,(i)}\bigr)\,\zeta_{j,o,U_1,U_2}^{tr,(i)},
\end{align*}
\ensuremath{k,l\in\{1,\ldots,Q\}.} The cost of approximating the family of the coefficient functions $a_{j,o,U_1,U_2}$, $j\in\left\{1,\ldots, J\right\}$, $(U_1,U_2)\in\mathcal{A}$, $\left|U_2\right|+\left|\mathcal{K}_2\right|+\frac{1}{2}\left|\mathcal{K}_1\right|\le 1$, $o\in\left\{1,2\right\}^{U_1}$, is of order
$O\bigl(Jm(m+1)NQ^{2}\bigr)$.

\subsection{Summary of the algorithm}
The algorithm consists of two phases:
training phase and testing phase.
In the training phase, we simulate
$N$ independent training paths $D_N^{tr}$
and construct regression estimates
$\hat a_{j,o,U_1,U_2}(\cdot,D_N^{tr})$
for the coefficients $a_{j,o,U_1,U_2}(\cdot)$.
In the testing phase,
independently from $D_N^{tr}$
we simulate $N_0$ independent testing paths
$(X_{\Delta,j\Delta}^{(i)})_{j=0,\ldots,J}$,
$i=1,\ldots,N_0$,
and build the Monte Carlo estimator
for $\EE[f(X_T)]$ as
\begin{equation}
\label{eq:0110a2}
\mathcal E=
\frac1{N_0}
\sum_{i=1}^{N_0}
\left(f(X^{(i)}_{\Delta,T})-\widehat M^{(2),trunc,(i)}_{\Delta,T}\right),
\end{equation}
where
\begin{align}\label{eq:19102016b1}
\widehat{M}^{(2),trunc,(i)}_{\Delta,T}\doteq\sum_{j=1}^J
\sum_{\substack{(U_1,U_2)\in\cA\\ \left|U_2\right|+\left|\mathcal{K}_2\right|+\frac{1}{2}\left|\mathcal{K}_1\right|\le 1}}
\sum_{o\in\{1,2\}^{U_1}}
\hat{a}_{j,o,U_1,U_2}(X^{(i)}_{\Delta,(j-1)\Delta},D_N^{tr})
\prod_{r\in U_1} H_{o_r}(\xi_j^{r,(i)})
\prod_{(k,l)\in U_2} V_j^{kl,(i)}
\end{align}
(cf.\ with~\eqref{eq:28042016a1}).
Due to the martingale transform structure
in~\eqref{eq:19102016b1}
(recall footnote~\ref{ft:19102016a1}
on page~\pageref{ft:19102016a1}),
we have
$\EE\left[\widehat M^{(2),trunc,(i)}_{\Delta,T}|D^{tr}_N\right]=0$,
hence
$\EE[\mathcal E|D^{tr}_N]
=\EE[f(X^{(i)}_{\Delta,T})-\widehat M^{(2),trunc,(i)}_{\Delta,T}|D^{tr}_N]
=\EE[f(X_{\Delta,T})]$,
and we obtain (cf.~\eqref{second:cv_var})
\begin{align*}
\Var[\mathcal E]
&=\EE[\Var(\mathcal E|D^{tr}_N)]
+\Var[\EE(\mathcal E|D^{tr}_N)]
=\EE[\Var(\mathcal E|D^{tr}_N)]\\
&=\frac1{N_0}
\EE\left[\Var\left(f(X^{(1)}_{\Delta,T})-\widehat M^{(2),trunc,(1)}_{\Delta,T}|D^{tr}_N\right)\right]
=\frac1{N_0}
\Var\left[f(X^{(1)}_{\Delta,T})-\widehat M^{(2),trunc,(1)}_{\Delta,T}\right].
\end{align*}
Summarising, we have
\begin{align}
\EE[\mathcal E]&=\EE[f(X_{\Delta,T})],
\notag\\
\Var[\mathcal E]&=\frac1{N_0}
\Var\left[f(X^{(1)}_{\Delta,T})-\widehat M^{(2),trunc,(1)}_{\Delta,T}\right].
\label{eq:19102016b3}
\end{align}
Notice that the result of~\eqref{eq:19102016b3}
indeed requires the computations above
and cannot be stated right from the outset
because the summands in~\eqref{eq:0110a2}
are dependent (through~$D^{tr}_N$).

This concludes the description
of the generic regression algorithm
for constructing the control variate.
Further details,
such as bounds for the right-hand side
of~\eqref{eq:19102016b3},
depend on a particular implementation,
i.e.\ on the quality of the chosen basis functions.

\section{Complexity analysis}
\label{sec:compl}

In this section we extend
the complexity analysis presented
in~\cite{belomestny2016variance}
to the case of the ``TRCV'' (truncated RCV) algorithm.
Below we only sketch the main results
for the second order schemes.
We make the following assumption (cf.~\cite{belomestny2016dual} and~\cite{belomestny2016stratified}):
\begin{itemize}
\item[(B5)]
The functions \(a_{j,o,U_1,U_2}\left(x\right)\) can be
well approximated by the functions from
$\Psi_{Q}\doteq\text{span}\left(\left\{\psi_{1},\ldots,\psi_{Q}\right\}\right)$,
in the sense that there are constants
$\kappa>0$ and $C_\kappa>0$ such that
\begin{align*}
\inf_{g\in\Psi_{Q}}\int_{\mathbb{R}^d}\left(a_{j,o,U_1,U_2}\left(x
\right)-g\left(x\right)\right)^2\,\PP_{\Delta,j-1}(dx)\leq  \frac{C_\kappa}{Q^{\kappa}},
\end{align*}
where $\PP_{\Delta,j-1}$
denotes the distribution of~$X_{\Delta,(j-1)\Delta}$.
\end{itemize}

\begin{remark}\label{rem:1301a1}
Note that~(B5) is a natural condition to be satisfied
for good choices of $\Psi_Q$.
For instance, under appropriate assumptions,
in the case of
\emph{piecewise polynomial regression}
as described in~\cite{belomestny2016variance}, 
(B5)~is satisfied with
$\kappa=\frac{2\nu(p+1)}{2d(p+1)+d\nu}$,
where the parameters $p$ and $\nu$
are explained in~\cite{belomestny2016variance}.
\end{remark}

In Lemma~\ref{th:2104a1}
below we present an $L^2$-upper bound
for the estimation error of the TRCV algorithm.
To this end, we need to describe more precisely,
how exactly the regression-based approximations
$\tilde a_{j,o,U_1,U_2}$ are constructed:

Let functions $\hat a_{j,o,U_1,U_2}(x)$
be obtained by regression onto the set of basis functions
$\left\{\psi_{1},\ldots,\psi_{Q}\right\}$,
while the approximations
$\tilde a_{j,o,U_1,U_2}(x)$ of the TRCV algorithm
be the truncated estimates, which are defined as follows
\begin{align}
\label{trunc_A}
\tilde a_{j,o,U_1,U_2}(x)\doteq
T_{\tilde A\Delta_{U_1,U_2}}\hat a_{j,o,U_1,U_2}(x)\doteq
\begin{cases}
\hat a_{j,o,U_1,U_2}(x)&\text{if }|\hat a_{j,o,U_1,U_2}(x)|\le\tilde A\Delta_{U_1,U_2},\\
\tilde A\Delta_{U_1,U_2}\sgn\hat a_{j,o,U_1,U_2}(x)&\text{otherwise}
\end{cases},
\end{align}
where $\Delta_{U_1,U_2}$ and $\tilde A$ are given in~\eqref{delta_u} and~\eqref{tilde_A}).

\begin{lemma}\label{th:2104a1}
Under~(B1)--(B5), we have
\begin{align}\label{eq:2104a2}
\E\|\tilde a_{j,o,U_1,U_2}-a_{j,o,U_1,U_2}\|^2_{L^2(\PP_{\Delta,j-1})}
\le\tilde{c}(\Sigma+\tilde A^2\Delta_{U_1,U_2}^2(\log N+1))\frac{Q}{N}
+\frac{8\,C_\kappa}{Q^\kappa},
\end{align}
where $\tilde{c}$ is a universal constant.
\end{lemma}
Notice that the expectation in the left-hand side of~\eqref{eq:2104a2}
means averaging over the randomness in~$D_N^{tr}$.

Let $(X_{\Delta,j\Delta})_{j=0,\ldots,J}$
be a testing path, which is independent of the
training paths $D_N^{tr}$.
We define
\begin{align}\label{eq:0101a1}
\widetilde{M}^{(2),trunc}_{\Delta,T}\doteq\sum_{j=1}^J
\sum_{\substack{(U_1,U_2)\in\cA\\ \left|U_2\right|+\left|\mathcal{K}_2\right|+\frac{1}{2}\left|\mathcal{K}_1\right|\le 1}}
\sum_{o\in\{1,2\}^{U_1}}
\tilde{a}_{j,o,U_1,U_2}(X_{\Delta,(j-1)\Delta},D_N^{tr})
\prod_{r\in U_1} H_{o_r}(\xi_j^{r})
\prod_{(k,l)\in U_2} V_j^{kl}
\end{align}
(cf.~\eqref{cv:second_trunc}).
Lemma~\ref{th:2104a1} now allows to bound the variance
$\Var[f(X_{\Delta,T})-\widetilde M_{\Delta,T}^{(2),trunc}]$
from above.\footnote{\label{ft:1301a1}Notice that the variance
of the TRCV estimate
$\frac{1}{N_0}\sum_{i=1}^{N_0}
\left[f\Bigl(X_{\Delta,T}^{(i)}\Bigr)-\tilde M_{\Delta,T}^{(2),trunc,(i)}\right]$
with $N_0$ testing paths
is $\frac1{N_0}\Var[f(X_{\Delta,T})-\tilde M_{\Delta,T}^{(2),trunc}]$ (cf.~\eqref{eq:19102016b3}).}

\begin{theorem}
\label{th:2106a1}
Under~(B1)--(B5), it holds
\begin{align*}
\Var[f(X_{\Delta,T})-
\widetilde{M}^{(2),trunc}_{\Delta,T}]
\lesssim \Delta^2+ Jm(m+1)\left(\tilde{c}(\Sigma+\tilde A^2\Delta(\log N+1))\frac{Q}{N}
+\frac{8\,C_\kappa}{Q^\kappa}\right).
\end{align*}
\end{theorem}

\subsection{Complexity of the TRCV approach}
\label{sec:trcv}
Let us study the complexity of the TRCV approach. The overall cost is of order $JQ\max\left\{NQ,N_0\right\}$, provided that we only track the constants which tend to $\infty$ when $\varepsilon\searrow 0$ with $\varepsilon$ being the accuracy to be achieved. That is, the constants, such as $d,m,\kappa,C_\kappa$, are ignored.
We have the following constraints
\begin{align}
\label{trcv_constr}
\max\left\{\frac{1}{J^4},\frac{1}{J^2 N_0},
\frac{JQ}{NN_0},
\frac{J}{Q^\kappa N_0}\right\}\lesssim\varepsilon^2,
\end{align}
where the first term comes from the squared bias of the estimator and the remaining three ones come from the variance of the estimator
(see Theorem~\ref{th:2106a1}
as well as footnote~\ref{ft:1301a1}
on page~\pageref{ft:1301a1}).
We get the following result.

\begin{theorem}
\label{trcv_compl}
For the TRCV approach with the second order weak schemes
under (B1)--(B5),
it is optimal to choose the orders of parameters as follows
(cf.~\cite{belomestny2016stratified})
\begin{gather*}
J\asymp\varepsilon^{-\frac{1}{2}},\quad Q\asymp \varepsilon^{-\frac{5}{4\kappa+4}},\quad N\asymp \varepsilon^{-\frac{5}{4}},\quad 
N_0\asymp NQ\asymp \varepsilon^{-\frac{5\kappa+10}{4\kappa+4}},
\end{gather*}
provided that $\kappa>1$.\footnote{\label{fn:trcv}Performing
the full complexity analysis via Lagrange multipliers one can see that these parameter values are
\emph{not} optimal if $\kappa\le 1$
(a Lagrange multiplier corresponding to
a ``$\le0$'' constraint is negative).
Recall that in the case of piecewise polynomial regression
(see~\cite{belomestny2016variance}
and recall Remark~\ref{rem:1301a1})
we have $\kappa=\frac{2\nu(p+1)}{2d(p+1)+d\nu}$.
Let us note that
in~\cite{belomestny2016variance} it is required
to choose the parameters $p$ and $\nu$ according to
$p>\frac{d-2}2$ and $\nu>\frac{2d(p+1)}{2(p+1)-d}$,
which implies that $\kappa>1$,
for $\kappa$ expressed via $p$ and $\nu$ by the above formula.}
Thus, we have for the complexity
\begin{align}
\label{eq:complt}
\mathcal{C}_{TRCV}&\asymp JNQ^2 \asymp JN_0Q\asymp \varepsilon^{-\frac{7\kappa+17}{4\kappa+4}}. 
\end{align}
\end{theorem}

\begin{remark}
(i) For the sake of comparison with the SMC and MLMC approaches,
we recall at this point that their complexities are
$$
\mathcal C_{SMC}\asymp\varepsilon^{-2.5}
\quad\text{and}\quad
\mathcal C_{MLMC}\asymp\varepsilon^{-2}
$$
at best (we are considering the second order scheme).

(ii) Complexity estimate~\eqref{eq:complt}
shows that one can go beyond the complexity order
$\varepsilon^{-2}$, provided that $\kappa>9$,
and that we can achieve the complexity order
$\varepsilon^{-1.75-\delta}$,
for arbitrarily small $\delta>0$,
provided $\kappa$ is large enough.

(iii) The complexity of the TRCV approach
is the same that we obtain for the RCV approach (where the ``complete'' control variate~\eqref{eq:28042016a1} is estimated), since the second constraint in~\eqref{trcv_constr}, which does not arise for the RCV approach, is the only inactive one in this case. That is why we truncated $M_{\Delta,T}^{(2),trunc}$ in~\eqref{cv:second_trunc} at the level $\left|U_2\right|+\left|\mathcal{K}_2\right|+\frac{1}{2}\left|\mathcal{K}_1\right|\le 1$.
For instance, if we had used a control variate of the form (cf.~\eqref{cv:first_trunc})
\begin{align*}
&\sum_{j=1}^J
\sum_{\substack{(U_1,U_2)\in\cA\\ \left|U_2\right|+\left|\mathcal{K}_2\right|+\frac{1}{2}\left|\mathcal{K}_1\right|= \frac{1}{2}}}
\sum_{o\in\{1,2\}^{U_1}}
a_{j,o,U_1,U_2}(X_{\Delta,(j-1)\Delta})
\prod_{r\in U_1} H_{o_r}(\xi_j^r)
\prod_{(k,l)\in U_2} V_j^{kl}\\
=&\sum_{j=1}^J\sum_{i=1}^ma_{j,1,i,\emptyset}(X_{\Delta,(j-1)\Delta})
\xi_j^i
\end{align*}
with $a_{j,1,i,\emptyset}(x)=\EE\left[f(X_{\Delta,T})\xi_j^i\left|\right.X_{\Delta,(j-1)\Delta}=x\right]$, 
the bound for the variance in~\eqref{second:cv_var}
would have been of order $\Delta$ and due to the resulting constraint $\frac{1}{JN_0}\lesssim\varepsilon^2$, we would have obtained worse complexities than $\varepsilon^{-2}$, since $\mathcal{C}_{TRCV}\gtrsim JN_0$.
\end{remark}

\section{Numerical results}
\label{sec:num}
The results below are based on program codes written and vectorised in MATLAB and running on a Linux 64-bit operating system.

Let us consider the following SDE for $d=m=5$ (cf.~\cite{belomestny2016variance})
\begin{align}
\notag
dX_t^i&=-\sin\left(X_t^i\right)\cos^3\left(X_t^i\right)dt+\cos^2\left(X_t^i\right)dW_t^i,\quad X_0^i=0,\quad i\in\left\{1,2,3,4\right\},\\
\label{5d_sde}
dX_t^5&=\sum_{i=1}^4\left[-\frac{1}{2}\sin\left(X_t^i\right)\cos^2\left(X_t^i\right)dt+\cos\left(X_t^i\right)
dW_t^i\right]+dW_t^5,\quad X_0^5=0.
\end{align}
The solution of~\eqref{5d_sde} is given by
\begin{align*}
X_t^i&=\arctan\left(W_t^i\right),\quad i\in\left\{1,2,3,4\right\},\\
X_t^5&=\sum_{i=1}^4\operatorname{arsinh}\left(W_t^i\right)
+W_t^5.
\end{align*}
for $t\in\left[0,1\right]$. Further, we consider the functional 
\begin{align*}
f(x)=\cos\left(\sum_{i=1}^5x^i\right)-20\sum_{i=1}^4\sin\left(x^i\right),
\end{align*}
that is, we have
\begin{align*}
\mathbb{E}\left[f\left(X_1\right)\right]&=\left(\mathbb{E}\left[\cos\left(
\arctan\left(W_1^1\right)+\operatorname{arsinh}\left(W_1^1\right)\right)\right]\right)^4
\mathbb{E}\left[\cos\left(W_1^5\right)\right]
\approx 0.002069.
\end{align*}

Here we consider weak schemes of the second order
and compare the numerical performances
of the SMC, MLMC, RCV, TRCV and TSRCV approaches.
The latter one is the truncated version
of the SRCV approach of~\cite{belomestny2016stratified}.
Like the RCV algorithm,
the SRCV one is based on~\eqref{eq:28042016a1},
the difference is only in how to implement
the approximations of the coefficients
$a_{j,o,U_1,U_2}$ in practice
(while the RCV algorithm is a direct Monte Carlo regression,
in the SRCV algorithm the regression is combined
with a kind of ``stratification''; see~\cite{belomestny2016stratified} for more detail).
Therefore, the idea of the truncation
(i.e.\ replacing~\eqref{eq:28042016a1} with~\eqref{cv:second_trunc})
applies also to the SRCV approach
and gives us the TSRCV one.

For simplicity we implemented a global regression for the RCV, TRCV and TSRCV approaches
(i.e.\ the one without considering the truncation operator in~\eqref{trunc_A}, as a part of the general description in Section~\ref{sec:regr}). 
More precisely, we use
quadratic polynomials (that is $\prod_{i=1}^5x_i^{l_i}$,
where $l_1,\ldots l_5\in\left\{0,1,2\right\}$
and $\sum_{l=1}^5l_i\leq 2$) as well as $f$ as basis functions, hence $\Psi_Q$ consists of $Q=\binom{7}{5}+1=22$ basis functions. \par
Note that we do not need to consider random variables $V_j^{kl}$
in the second order weak scheme,
since $\cL^k\sigma^{rl}(x)=0$ for $k\ne l$ (see~\eqref{eq:2002a6}).
This gives us less terms for the RCV approach, namely $3^m-1=242$ rather than $3^m2^\frac{m(m-1)}{2}-1=248831$ terms in~\eqref{eq:28042016a1} (the factor $2^{\frac{m(m-1)}2}\equiv 1024$ is no longer present). As for the TRCV
and TSRCV approaches,
this gives us only $\frac{m(m+3)}{2}=20$ compared to $m(m+1)=30$ terms in~\eqref{cv:second_trunc}.

\begin{figure}[htb!]
\includegraphics[width=0.70\textwidth]{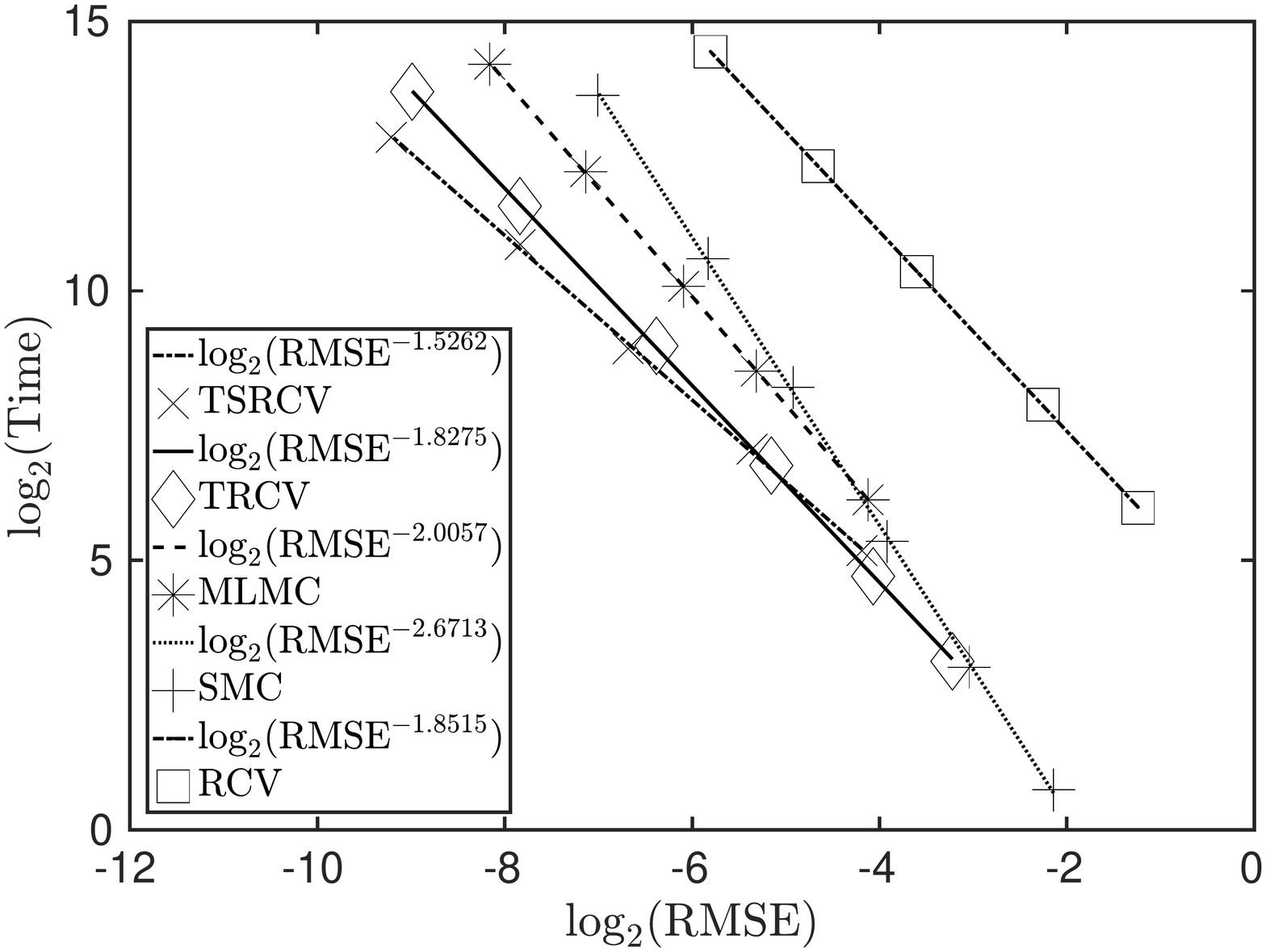}
\caption{Numerical complexities of the RCV, TRCV, TSRCV, SMC and MLMC approaches.}
\label{compld}
\end{figure}

We choose $\kappa=1.2$, which is related to the piecewise polynomial regression with polynomial degree $p=2$ (comparable to our setting) and the limiting case $\nu\to\infty$ (see footnote~\ref{fn:trcv} on page~\pageref{fn:trcv}). Moreover, for each $\varepsilon=2^{-i}$, $i\in\left\{2,3,4,5,6\right\}$,
we set the parameters $J$, $N$ and $N_0$ for the RCV,
TRCV and TSRCV approaches as follows
(compare with the formulas in Subsection~\ref{sec:trcv}):
\begin{align*}
J=\left\lceil \varepsilon^{-0.5}\right\rceil,\qquad
N=\begin{cases}
512\cdot\lceil \varepsilon^{-1.25}\rceil,&\text{RCV, TRCV},\\
2048\cdot\lceil \varepsilon^{-1.25}\rceil,&\text{TSRCV},
\end{cases}
\qquad 
N_0=512\cdot\lceil \varepsilon^{-1.82}\rceil.
\end{align*}
The factors 512 and 2048 are here for stability purposes.
For the TRCV and SMC algorithms we additionally consider
$\varepsilon=2^{-7}$,
which produces a picture with approximately equal
maximal computational time
(that is, the time corresponding to the best accuracy)
for all algorithms.
Next we estimate the numerical complexity for the RCV, TRCV and TSRCV approaches by means of 100 independent simulations and compare it with the one for the SMC and MLMC approach, for which we use the same output as in~\cite{belomestny2016variance}. As can be seen from Figure~\ref{compld}, the estimated numerical complexity is about $\text{RMSE}^{-1.85}$ for the RCV approach, $\text{RMSE}^{-1.83}$ for the TRCV approach,
$\text{RMSE}^{-1.53}$ for the TSRCV approach,
$\text{RMSE}^{-2.67}$ for the SMC approach and $\text{RMSE}^{-2.01}$ for the MLMC approach,
which we get by regressing the log-time (logarithmic computing time of the whole algorithm in seconds) vs.\ log-RMSE. Beyond the numerical complexities we observe that the truncation effect from RCV algorithm to its truncated versions is huge. While we have poor results for the RCV approach (as in~\cite{belomestny2016variance}), i.e. in this region of $\varepsilon$-values the RCV approach is numerically outperformed by the other ones, the TRCV and TSRCV approaches work best
(even better than the SMC and MLMC approaches).

\section{Proofs}
\label{sec:proofs}
\subsection*{Proof of Theorem~\protect\ref{smooth:func_assump_general}}
We begin with the following remark.
Assumptions~\eqref{smooth:Phi_assump_1}
and~\eqref{smooth:Phi_assump_2 new}
together with the Cauchy-Schwarz inequality
$|\EE[XY|\cG]|\le\sqrt{\EE[X^2|\cG]\,\EE[Y^2|\cG]}$
imply that the following generalisation
of~\eqref{smooth:Phi_assump_2 new} is satisfied:
for any $n_1,n_2\in\NN$, $\alpha,\beta\in\NN_0^d$,
with
$1\le\left|\alpha\right|\le K$,
$1\le\left|\beta\right|\le K$,
$\alpha\ne\beta$,
it holds
\begin{align}
\label{smooth:Phi_assump_2}
\left|\EE\left[\left.\left(D^{\alpha}\Phi{}_{\Delta,l+1}^{k}(G_{l,j}(x))\right)^{n_1}\left(D^{\beta}\Phi{}_{\Delta,l+1}^{k}(G_{l,j}(x))\right)^{n_2}\right|\mathcal{G}_{l}\right]\right|\leq C_{n_1,n_2}\Delta
\end{align}
for some appropriate constants $C_{n_1,n_2}>0$.

Let us 
begin with the case $K=1$. We have for some $k,r\in\left\{1,\ldots,d\right\}$
$$
\frac{\partial}{\partial x^r}G_{l+1,j}^k(x)=\sum_{s=1}^d\frac{\partial}{\partial x^s}\Phi_{\Delta,l+1}^{k}(G_{l,j}(x))\frac{\partial}{\partial x^r}G_{l,j}^s(x)=:\sum_{s=1}^d\gamma_s
$$
and $\frac{\partial}{\partial x^r}G_{j+1,j}^s(x)=\frac{\partial}{\partial x^r}\Phi_{\Delta}^s\left(x,\xi_{j+1}\right)$, where $G_{l+1,j}^s$ and $\Phi_\Delta^s$, $s\in\left\{1,\ldots,d\right\}$, denote the $s$-th component of the functions $G_{l+1,j}$ and $\Phi_\Delta$.
Hence 
\begin{eqnarray*}
&&\EE\left[\left(\frac{\partial}{\partial x^r}G_{l+1,j}^k(x)\right)^{2}\right]\leq\EE\left[\gamma_k^2+\sum_{s:\,s\neq k}\left(2\gamma_k\gamma_s+(d-1)\gamma_s^2\right)\right].
\end{eqnarray*}
For an arbitrary $j\in\{0,\ldots,J-1\}$, denote 
\[
\rho_{l+1,n,1}^{r,s}:=\EE\left[\left(\frac{\partial}{\partial x^r}G_{l+1,j}^s(x)\right)^{n}\right],
\]
then, due to~\eqref{smooth:Phi_assump_1} and~\eqref{smooth:Phi_assump_2}, we get for $l=j,\ldots,J-1$,
\begin{eqnarray*}
\rho_{l+1,2,1}^{r,k}  \leq  (1+A_{2}\Delta)\rho_{l,2,1}^{r,k}
+\sum_{s:\,s\neq k}\left(C_{1,1}\Delta(\rho_{l,2,1}^{r,k}+\rho_{l,2,1}^{r,s})+(d-1)B_{2}\Delta\rho_{l,2,1}^{r,s}\right).
\end{eqnarray*}
Further, denote 
\[
\rho_{l+1,n,1}^{r}:=\sum_{s=1}^d\rho_{l+1,n,1}^{r,s},
\]
then we get 
\begin{eqnarray*}
\rho_{l+1,2,1}^{r}  \leq  (1+A_{2}\Delta)\rho_{l,2,1}^{r}
+2(d-1)C_{1,1}\Delta\rho_{l,2,1}^{r}+(d-1)^2B_{2}\Delta\rho_{l,2,1}^{r}.
\end{eqnarray*}
This gives us
\[
\rho_{l+1,2,1}^{r}\leq (1+\kappa_{1}\Delta)\rho_{l,2,1}^{r}
\]
for some constant $\kappa_{1}>0$,
leading to
\begin{equation}\label{eq:13052017a1}
\rho_{l,2,1}^{r}\leq(1+\kappa_{1}\Delta)^{l-j-1}\rho_{j+1,2,1}^{r}
,\quad l=j+1,\ldots,J-1,
\end{equation}
where 
\[
\rho_{j+1,2,1}^{r}=\sum_{s=1}^d\EE\left[\left(\frac{\partial}{\partial x^r}\Phi_{\Delta}^s\left(x,\xi_{j+1}\right)\right)^{2}\right],
\]
which is bounded due to~\eqref{smooth:Phi_assump_1}.
Together with~\eqref{eq:13052017a1}
we obtain the boundedness of
$\{\rho^r_{J,2,1}:J\in\NN\}$
and hence the boundedness of
\begin{eqnarray*}
\left|\frac{\partial}{\partial x^r}q_j(x)\right|&\le&\sum_{s=1}^d\EE\left|\frac{\partial }{\partial x^s}f(G_{J,j}(x))\frac{\partial}{\partial x^r}G_{J,j}^s(x)\right|\\
&\le &\sum_{s=1}^d\sqrt{\EE\left[\left(\frac{\partial }{\partial x^s}f(G_{J,j}(x))\right)^2\right]\rho_{J,2,1}^{r,s}}\\
&\le &\sqrt{d\sum_{s=1}^d
\EE\left[\left(\frac{\partial }{\partial x^s}f(G_{J,j}(x))\right)^2\right]\rho_{J,2,1}^{r,s}}
\le\const\sqrt{\rho^r_{J,2,1}}
\end{eqnarray*}
for all $r\in\left\{1,\ldots,d\right\}$, since $f$ is assumed to be continuously differentiable with bounded partial derivatives.

Let us proceed with the case $K=2$. We have, due to $(\sum_{k=1}^da_k)^n\le d^{n-1}\sum_{k=1}^d a_k^n$,
\begin{eqnarray*}
&&\EE\left[\left(\frac{\partial}{\partial x^r}G_{l+1,j}^k(x)\right)^{4}\right]\\
&\leq&\EE\left[\gamma_k^4+\sum_{s:\,s\neq k}\left(4\gamma_k^3\gamma_s+6(d-1)\gamma_k^2\gamma_s^2+4(d-1)^2\gamma_k\gamma_s^3+(d-1)^3\gamma_s^4\right)\right]
\end{eqnarray*}
and thus, due to $4a^3b\le 3a^4+b^4$ and $2a^2b^2\le a^4+b^4$,
\begin{eqnarray*}
\rho_{l+1,4,1}^{r,k}
&\le & (1+A_4\Delta)\rho_{l,4,1}^{r,k}
+\sum_{s:\,s\neq k}\left(C_{3,1}\Delta(3\rho_{l,4,1}^{r,k}+\rho_{l,4,1}^{r,s})+3(d-1)C_{2,2}\Delta(\rho_{l,4,1}^{r,k}+\rho_{l,4,1}^{r,s})\right.\\
&&\phantom{(1+A_{4}\Delta)\rho_{l,4,1}^{r,k}
+\sum_{s:\,s\neq k}\left(\right.}\left.+(d-1)^2C_{1,3}\Delta(\rho_{l,4,1}^{r,k}+3\rho_{l,4,1}^{r,s})+ (d-1)^3B_{4}\Delta\rho_{l,4,1}^{r,s}\right).
\end{eqnarray*}
This gives us
\begin{eqnarray*}
\rho_{l+1,4,1}^{r}
&\le & (1+A_{4}\Delta)\rho_{l,4,1}^{r}
+4(d-1)C_{3,1}\Delta\rho_{l,4,1}^{r}+6(d-1)^2C_{2,2}\Delta\rho_{l,4,1}^{r}\\
&&+4(d-1)^3C_{1,3}\Delta\rho_{l,4,1}^{r}+ (d-1)^4B_{4}\Delta\rho_{l,4,1}^{r}.
\end{eqnarray*}
Hence, we obtain
\[
\rho_{l+1,4,1}^{r}\leq (1+\kappa_{2}\Delta)\rho_{l,4,1}^{r},
\]
for some constant $\kappa_{2}>0$, leading to 
\[
\rho_{l,4,1}^{r}\leq(1+\kappa_{2}\Delta)^{l-j-1}\rho_{j+1,4,1}^{r}
,\quad l=j+1,\ldots,J-1,
\]
where 
\[
\rho_{j+1,4,1}^{r}=\sum_{s=1}^d\EE\left[\left(\frac{\partial}{\partial x^r}\Phi_{\Delta}^s\left(x,\xi_{j+1}\right)\right)^{4}\right].
\]
Next, we have for some $k,o,r\in\left\{1,\ldots,d\right\}$
\begin{align*}
\frac{\partial^2}{\partial x^r\partial x^o}G_{l+1,j}^k(x)&=\sum_{s=1}^d\frac{\partial}{\partial x^s}\Phi_{\Delta,l+1}^k(G_{l,j}(x))\frac{\partial^2}{\partial x^r\partial x^o}G_{l,j}^s(x)\\
&\phantom{=}+\sum_{s,u=1}^d\frac{\partial^2}{\partial x^s\partial x^u}\Phi_{\Delta,l+1}^k(G_{l,j}(x))\frac{\partial}{\partial x^r}G_{l,j}^s(x)\frac{\partial}{\partial x^o}G_{l,j}^u(x)\\
&=:\sum_{s=1}^d\eta_{1,s}+\sum_{s,u=1}^d\eta_{2,s,u}
\end{align*}
and $\frac{\partial^2}{\partial x^r\partial x^o}G_{j+1,j}^s(x)=\frac{\partial^2}{\partial x^r\partial x^o}\Phi_{\Delta}^s\left(x,\xi_{j+1}\right)$.
Hence 
\begin{eqnarray*}
&&\EE\left[\left(\frac{\partial^2}{\partial x^r\partial x^o}G_{l+1,j}^k(x)\right)^{2}\right]\\
&\leq&\EE\left[\eta_{1,k}^2+\sum_{s:\,s\neq k}\left(2\eta_{1,k}\eta_{1,s}+(d-1)\eta_{1,s}^2\right)+2\sum_{s,u,v=1}^d\eta_{1,v}\eta_{2,s,u}+d^2\sum_{s,u=1}^d\eta_{2,s,u}^2\right]
\end{eqnarray*}
Denote 
\[
\rho_{l+1,n,2}^{r,o,s}=\EE\left[\left(\frac{\partial^2}{\partial x^r\partial x^o}G_{l+1,j}^s(x)\right)^{n}\right],
\]
then we get, due to 
\begin{align*}
2\EE\left[XYZ\right]&\le 2\sqrt{\EE\left[X^2\right]}\sqrt[4]{\EE\left[Y^4\right]}\sqrt[4]{\EE\left[Z^4\right]}\le \EE\left[X^2\right]+\sqrt{\EE\left[Y^4\right]}\sqrt{\EE\left[Z^4\right]}\\
&\le \EE\left[X^2\right]+\frac{1}{2}\left(\EE\left[Y^4\right]+\EE\left[Z^4\right]\right),
\end{align*}
as well as~\eqref{smooth:Phi_assump_1} and~\eqref{smooth:Phi_assump_2},
\begin{eqnarray*}
\rho_{l+1,2,2}^{r,o,k} & \leq & (1+A_{2}\Delta)\rho_{l,2,2}^{r,o,k}
+\sum_{s:\,s\neq k}\left(C_{1,1}\Delta(\rho_{l,2,2}^{r,o,k}+\rho_{l,2,2}^{r,o,s})+(d-1)B_{2}\Delta\rho_{l,2,2}^{r,o,s}\right)\\
&&+\sum_{s,u,v=1}^d C_{1,1}\Delta\left(\rho_{l,2,2}^{r,o,v}+\frac{1}{2}\left(\rho_{l,4,1}^{r,s}+\rho_{l,4,1}^{o,u}\right)\right)\\
&&+d^2\sum_{s,u=1}^dB_{2}\Delta\frac{1}{2}\left(\rho_{l,4,1}^{r,s}+\rho_{l,4,1}^{o,u}\right).
\end{eqnarray*}
Further, denote 
\[
\rho_{l+1,n,2}^{r,o}=\sum_{s=1}^d\rho_{l+1,n,2}^{r,o,s},
\]
then we get for $l=j+1,\ldots,J-1$,
\begin{eqnarray*}
\rho_{l+1,2,2}^{r,o} & \leq & (1+A_{2}\Delta)\rho_{l,2,2}^{r,o}
+2(d-1)C_{1,1}\Delta\rho_{l,2,2}^{r,o}+(d-1)^2B_{2}\Delta\rho_{l,2,2}^{r,o}\\
&&+d^3C_{1,1}\Delta\left(\rho_{l,2,2}^{r,o}+\frac{1}{2}\left(\rho_{l,4,1}^{r}+\rho_{l,4,1}^{o}\right)\right)+d^4B_2\Delta\frac{1}{2}\left(\rho_{l,4,1}^{r}+\rho_{l,4,1}^{o}\right).
\end{eqnarray*}
This gives us
\[
\rho_{l+1,2,2}^{r,o}\leq (1+\kappa_{3}\Delta)\rho_{l,2,2}^{r,o}+\kappa_4\Delta,
\]
for some constants $\kappa_{3},\kappa_4>0$,
leading to
\[
\rho_{l,2,2}^{r,o}\leq(1+\kappa_{3}\Delta)^{l-j-1}\rho_{j+1,2,2}^{r,o}+\kappa_5
,\quad l=j+1,\ldots,J-1,
\]
where $\kappa_5>0$ and
\[
\rho_{j+1,2,2}^{r,o}=\sum_{s=1}^d\EE\left[\left(\frac{\partial^2}{\partial x^r\partial x^o}\Phi_{\Delta}^s\left(x,\xi_{j+1}\right)\right)^{2}\right].
\]
Thus, we obtain the boundedness of 
\begin{eqnarray*}
\left|\frac{\partial^2}{\partial x^r \partial x^o}q_j(x)\right|&\le&\sum_{s=1}^d\EE\left|\frac{\partial }{\partial x^s}f(G_{J,j}(x))\frac{\partial^2}{\partial x^r\partial x^o}G_{J,j}^s(x)\right|\\
&&+\sum_{s,u=1}^d\EE\left|\frac{\partial^2 }{\partial x^s\partial x^u}f(G_{J,j}(x))\frac{\partial}{\partial x^r}G_{J,j}^s(x)\frac{\partial}{\partial x^o}G_{J,j}^u(x)\right|\\
&\le &\sum_{s=1}^d\sqrt{\EE\left[\left(\frac{\partial }{\partial x^s}f(G_{J,j}(x))\right)^2\right]\rho_{J,2,2}^{r,o,s}}\\
&&+\sum_{s,u=1}^d\sqrt{\EE\left[\left(\frac{\partial^2 }{\partial x^s\partial x^u}f(G_{J,j}(x))\right)^2\right]}\sqrt[4]{\rho_{J,4,1}^{r,s}\rho_{J,4,1}^{o,u}}
\end{eqnarray*}
for all $r,o\in\left\{1,\ldots,d\right\}$, since $f$ is assumed to be twice continuously differentiable with bounded partial derivatives up to order $2$.

Let us proceed with the final case $K=3$.
We have
\begin{eqnarray*}
&&\EE\left[\left(\frac{\partial}{\partial x^r}G_{l+1,j}^k(x)\right)^{6}\right]\\
&\leq&\EE\left[\gamma_k^6+\sum_{s:\,s\neq k}\left(6\gamma_k^5\gamma_s+15(d-1)\gamma_k^4\gamma_s^2+20(d-1)^2\gamma_k^3\gamma_s^3+15(d-1)^3\gamma_k^2\gamma_s^4\right.\right.\\
&&\phantom{\EE\left[\gamma_k^6+\sum_{s:\,s\neq k}\left(\right.\right.}\left.\left.+6(d-1)^4\gamma_k\gamma_s^5+(d-1)^5\gamma_s^6\right)\right]
\end{eqnarray*}
and thus, due to $6a^5b\le 5a^6+b^6$, $3a^4b^2\le 2a^6+b^6$ and $2a^3b^3\le a^6+b^6$,
\begin{eqnarray*}
\rho_{l+1,6,1}^{r,k}&\le & (1+A_{6}\Delta)\rho_{l,6,1}^{r,k}\\
&&+\sum_{s:\,s\neq k}\left(C_{5,1}\Delta(5\rho_{l,6,1}^{r,k}+\rho_{l,6,1}^{r,s})+5(d-1)C_{4,2}\Delta(2\rho_{l,6,1}^{r,k}+\rho_{l,6,1}^{r,s})\right.\\
&&\phantom{+\sum_{s:\,s\neq k}\left(\right.}+10(d-1)^2C_{3,3}\Delta(\rho_{l,6,1}^{r,k}+\rho_{l,6,1}^{r,s})+5(d-1)^3C_{2,4}\Delta(\rho_{l,6,1}^{r,k}+2\rho_{l,6,1}^{r,s})\\
&&\phantom{+\sum_{s:\,s\neq k}\left(\right.}\left.+(d-1)^4C_{1,5}\Delta(\rho_{l,6,1}^{r,k}+5\rho_{l,6,1}^{r,s})+ (d-1)^5B_{6}\Delta\rho_{l,6,1}^{r,s}\right).
\end{eqnarray*}
This gives us
\begin{eqnarray*}
\rho_{l+1,6,1}^{r}&\le & (1+A_{6}\Delta)\rho_{l,6,1}^{r}
+6(d-1)C_{5,1}\Delta\rho_{l,6,1}^{r}+15(d-1)^2C_{4,2}\Delta\rho_{l,6,1}^{r}\\
&&+20(d-1)^3C_{3,3}\Delta\rho_{l,6,1}^{r}+15(d-1)^4C_{2,4}\Delta\rho_{l,6,1}^{r}+6(d-1)^5C_{1,5}\Delta\rho_{l,6,1}^{r}\\
&&+ (d-1)^6B_{6}\Delta\rho_{l,6,1}^{r}.
\end{eqnarray*}
Hence, we obtain
\[
\rho_{l+1,6,1}^{r}\leq (1+\kappa_{6}\Delta)\rho_{l,6,1}^{r}
\]
for some constant $\kappa_{6}>0$, leading to 
\[
\rho_{l,6,1}^{r}\leq(1+\kappa_{6}\Delta)^{l-j-1}\rho_{j+1,6,1}^{r}
,\quad l=j+1,\ldots,J-1,
\]
where 
\[
\rho_{j+1,6,1}^{r}=\sum_{s=1}^d\EE\left[\left(\frac{\partial}{\partial x^r}\Phi_{\Delta}^s\left(x,\xi_{j+1}\right)\right)^{6}\right]
\]
Moreover, we have
\begin{eqnarray*}
&&\EE\left[\left(\frac{\partial}{\partial x^r}G_{l+1,j}^k(x)\right)^{8}\right]\\
&\leq&\EE\left[\gamma_k^8+\sum_{s:\,s\neq k}\left(8\gamma_k^7\gamma_s+28(d-1)\gamma_k^6\gamma_s^2+56(d-1)^2\gamma_k^5\gamma_s^3+70(d-1)^3\gamma_k^4\gamma_s^4\right.\right.\\
&&\phantom{\EE\left[\right.\gamma_k^8+\sum_{s:\,s\neq k}\left(\right.}\left.\left.+56(d-1)^4\gamma_k^3\gamma_s^5+28(d-1)^5\gamma_k^2\gamma_s^6+8(d-1)^6\gamma_k\gamma_s^7+(d-1)^7\gamma_s^8\right)\right]
\end{eqnarray*}
and thus, due to $8a^7b\le 7a^8+b^8$, $4a^6b^2\le 3a^8+b^8$, $8a^5b^3\le 5a^8+3b^8$ and $2a^4b^4\le a^8+b^8$,
\begin{eqnarray*}
\rho_{l+1,8,1}^{r,k}
&\le & (1+A_{8}\Delta)\rho_{l,8,1}^{r,k}\\
&&
+\sum_{s:\,s\neq k}\left(C_{7,1}\Delta(7\rho_{l,8,1}^{r,k}+\rho_{l,8,1}^{r,s})+7(d-1)C_{6,2}\Delta(3\rho_{l,8,1}^{r,k}+\rho_{l,8,1}^{r,s})\right.\\
&&\phantom{+\sum_{s:\,s\neq k}\left(\right.}+7(d-1)^2C_{5,3}\Delta(5\rho_{l,8,1}^{r,k}+3\rho_{l,8,1}^{r,s})+35(d-1)^3C_{4,4}\Delta(\rho_{l,8,1}^{r,k}+\rho_{l,8,1}^{r,s})\\
&&\phantom{+\sum_{s:\,s\neq k}\left(\right.}+7(d-1)^4C_{3,5}\Delta(3\rho_{l,8,1}^{r,k}+5\rho_{l,8,1}^{r,s})+7(d-1)^5C_{2,6}\Delta(\rho_{l,8,1}^{r,k}+3\rho_{l,8,1}^{r,s})\\
&&\phantom{+\sum_{s:\,s\neq k}\left(\right.}\left.+(d-1)^6C_{1,7}\Delta(\rho_{l,8,1}^{r,k}+7\rho_{l,8,1}^{r,s})+ (d-1)^7B_{8}\Delta\rho_{l,8,1}^{r,s}\right).
\end{eqnarray*}
This gives us
\begin{eqnarray*}
\rho_{l+1,8,1}^{r}
&\le & (1+A_{8}\Delta)\rho_{l,8,1}^{r}
+8(d-1)C_{7,1}\Delta\rho_{l,8,1}^{r}+28(d-1)^2C_{6,2}\Delta\rho_{l,8,1}^{r}\\
&&+56(d-1)^3C_{5,3}\Delta\rho_{l,8,1}^{r}+70(d-1)^4C_{4,4}\Delta\rho_{l,8,1}^{r}\\
&&+56(d-1)^5C_{3,5}\Delta\rho_{l,8,1}^{r}+28(d-1)^6C_{2,6}\Delta\rho_{l,8,1}^{r}\\
&&+8(d-1)^7C_{1,7}\Delta\rho_{l,8,1}^{r}+ (d-1)^8B_{8}\Delta\rho_{l,8,1}^{r}.
\end{eqnarray*}
Hence, we obtain
\[
\rho_{l+1,8,1}^{r}\leq (1+\kappa_{7}\Delta)\rho_{l,8,1}^{r},
\]
for some constant $\kappa_{7}>0$, leading to 
\[
\rho_{l,8,1}^{r}\leq(1+\kappa_{7}\Delta)^{l-j-1}\rho_{j+1,8,1}^{r}
,\quad l=j+1,\ldots,J-1,
\]
where 
\[
\rho_{j+1,8,1}^{r}=\sum_{s=1}^d\EE\left[\left(\frac{\partial}{\partial x^r}\Phi_{\Delta}^s\left(x,\xi_{j+1}\right)\right)^{8}\right].
\]
Moreover, we have
\begin{eqnarray*}
&&\EE\left[\left(\frac{\partial^2}{\partial x^r\partial x^o}G_{l+1,j}^k(x)\right)^{4}\right]\\
&\leq&\EE\left[\eta_{1,k}^4+\sum_{s:\,s\neq k}\left(4\eta_{1,k}^3\eta_{1,s}+6(d-1)\eta_{1,k}^2\eta_{1,s}^2+4(d-1)^2\eta_{1,k}\eta_{1,s}^3+(d-1)^3\eta_{1,s}^4\right)\right.\\
&&\phantom{\EE\left[\right.}\left.+\sum_{s,u,v=1}^d\left(4d^2\eta_{1,v}^3\eta_{2,s,u}+6d^3\eta_{1,v}^2\eta_{2,s,u}^2+4d^4\eta_{1,v}\eta_{2,s,u}^3\right)+d^6\sum_{s,u=1}^d\eta_{2,s,u}^4\right]
\end{eqnarray*}
and thus, due to $4a^3bc\le3a^4+\frac{1}{2}\left(b^8+c^8\right)$, $2a^2b^2c^2\le a^4+\frac{1}{2}\left(b^8+c^8\right)$ and $4ab^3c^3\le a^4+\frac{3}{2}\left(b^8+c^8\right)$,
\begin{eqnarray*}
\rho_{l+1,4,2}^{r,o,k}&\le & (1+A_{4}\Delta)\rho_{l,4,2}^{r,o,k}\\
&&
+\sum_{s:\,s\neq k}\left(C_{3,1}\Delta(3\rho_{l,4,2}^{r,o,k}+\rho_{l,4,2}^{r,o,s})+3(d-1)C_{2,2}\Delta(\rho_{l,4,2}^{r,o,k}+\rho_{l,4,2}^{r,o,s})\right.\\
&&\phantom{+\sum_{s:\,s\neq k}\left(\right.}\left.+(d-1)^2C_{1,3}\Delta(\rho_{l,4,2}^{r,o,k}+3\rho_{l,4,2}^{r,o,s})+ (d-1)^3B_{4}\Delta\rho_{l,4,2}^{r,o,s}\right)\\
&&+\sum_{s,u,v=1}^d\left(d^2C_{3,1}\Delta\left(3\rho_{l,4,2}^{r,o,v}+\frac{1}{2}\left(\rho_{l,8,1}^{r,s}+\rho_{l,8,1}^{o,u}\right)\right)\right.\\
&&\phantom{+\sum_{s,u,v=1}^d\left(\right.}+3d^3C_{2,2}\Delta\left(\rho_{l,4,2}^{r,o,v}+\frac{1}{2}\left(\rho_{l,8,1}^{r,s}+\rho_{l,8,1}^{o,u}\right)\right)\\
&&\phantom{+\sum_{s,u,v=1}^d\left(\right.}\left.+d^4C_{1,3}\Delta\left(\rho_{l,4,2}^{r,o,v}+\frac{3}{2}\left(\rho_{l,8,1}^{r,s}+\rho_{l,8,1}^{o,u}\right)\right)\right)\\
&&+d^6\sum_{s,u=1}^dB_{4}\Delta\frac{1}{2}\left(\rho_{l,8,1}^{r,s}+\rho_{l,8,1}^{o,u}\right).
\end{eqnarray*}
This gives us
\begin{eqnarray*}
\rho_{l+1,4,2}^{r,o}
&\le & (1+A_{4}\Delta)\rho_{l,4,2}^{r,o}
+4(d-1)C_{3,1}\Delta\rho_{l,4,2}^{r,o}+6(d-1)^2C_{2,2}\Delta\rho_{l,4,2}^{r,o}\\
&&+4(d-1)^3C_{1,3}\Delta\rho_{l,4,1}^{r}+ (d-1)^4B_{4}\Delta\rho_{l,4,2}^{r,o}\\
&&+d^5C_{3,1}\Delta\left(3\rho_{l,4,2}^{r,o}+\frac{1}{2}\left(\rho_{l,8,1}^{r}+\rho_{l,8,1}^{o}\right)\right)\\
&&+3d^6C_{2,2}\Delta\left(\rho_{l,4,2}^{r,o}+\frac{1}{2}\left(\rho_{l,8,1}^{r}+\rho_{l,8,1}^{o}\right)\right)\\
&&+d^7C_{1,3}\Delta\left(\rho_{l,4,2}^{r,o}+\frac{3}{2}\left(\rho_{l,8,1}^{r}+\rho_{l,8,1}^{o}\right)\right)+d^8B_{4}\Delta\frac{1}{2}\left(\rho_{l,8,1}^{r}+\rho_{l,8,1}^{o}\right).
\end{eqnarray*}
Hence, we obtain
\[
\rho_{l+1,4,2}^{r,o}\leq (1+\kappa_{8}\Delta)\rho_{l,4,2}^{r,o}+\kappa_9\Delta,
\]
for some constants $\kappa_{8},\kappa_9>0$, leading to 
\[
\rho_{l,4,2}^{r,o}\leq(1+\kappa_{8}\Delta)^{l-j-1}\rho_{j+1,4,2}^{r,o}+\kappa_{10}
,\quad l=j+1,\ldots,J-1,
\]
where $\kappa_{10}>0$ and
\[
\rho_{j+1,4,2}^{r,o}=\sum_{s=1}^d\EE\left[\left(\frac{\partial^2}{\partial x^r\partial x^o}\Phi_{\Delta}^s\left(x,\xi_{j+1}\right)\right)^{4}\right].
\]
Next, we have for some $k,o,r,z\in\left\{1,\ldots,d\right\}$
\begin{align*}
\phantom{=:}&\frac{\partial^3}{\partial x^r\partial x^o\partial x^z}G_{l+1,j}^k(x)\\
=\phantom{:}&\sum_{s=1}^d\frac{\partial}{\partial x^s}\Phi_{\Delta,l+1}^k(G_{l,j}(x))\frac{\partial^3}{\partial x^r\partial x^o\partial x^z}G_{l,j}^s(x)\\
&+\sum_{s,u=1}^d\frac{\partial^2}{\partial x^s\partial x^u}\Phi_{\Delta,l+1}^k(G_{l,j}(x))\left(\frac{\partial^2}{\partial x^r\partial x^o}G_{l,j}^s(x)\frac{\partial}{\partial x^z}G_{l,j}^u(x)+\frac{\partial^2}{\partial x^r\partial x^z}G_{l,j}^s(x)\frac{\partial}{\partial x^o}G_{l,j}^u(x)\right.\\
&\phantom{+\sum_{s,u=1}^d\frac{\partial^2}{\partial x^s\partial x^u}\Phi_{\Delta,l+1}^k(G_{l,j}(x))\left(\right.}\left.+\frac{\partial}{\partial x^r}G_{l,j}^s(x)\frac{\partial^2}{\partial x^o\partial x^z}G_{l,j}^u(x)\right)\\
&+\sum_{s,u,v=1}^d\frac{\partial^3}{\partial x^s\partial x^u\partial x^v}\Phi_{\Delta,l+1}^k(G_{l,j}(x))\frac{\partial}{\partial x^r}G_{l,j}^s(x)\frac{\partial}{\partial x^o}G_{l,j}^u(x)\frac{\partial}{\partial x^z}G_{l,j}^v(x)\\
=:&\sum_{s=1}^d\psi_{1,s}+\sum_{s,u=1}^d\psi_{2,s,u}+\sum_{s,u,v=1}^d\psi_{3,s,u,v}
\end{align*}
and $\frac{\partial^3}{\partial x^r\partial x^o\partial x^z}G_{j+1,j}^s(x)=\frac{\partial^3}{\partial x^r\partial x^o\partial x^z}\Phi_{\Delta}^s\left(x,\xi_{j+1}\right)$.
Hence 
\begin{eqnarray*}
&&\EE\left[\left(\frac{\partial^3}{\partial x^r\partial x^o\partial x^z}G_{l+1,j}^k(x)\right)^{2}\right]\\
&\leq&\EE\left[\psi_{1,k}^2+\sum_{s:\,s\neq k}\left(2\psi_{1,k}\psi_{1,s}+(d-1)\psi_{1,s}^2\right)+2\sum_{s,u,v=1}^d\psi_{1,v}\psi_{2,s,u}\right.\\
&&\phantom{\EE\left[\right.}\left.+2\sum_{s,u,v,w=1}^d\psi_{1,w}\psi_{3,s,u,v}+2d^2\sum_{s,u=1}^d\psi_{2,s,u}^2+2d^3\sum_{s,u,v=1}^d\psi_{3,s,u,v}^2\right]
\end{eqnarray*}
Denote 
\[
\rho_{l+1,n,3}^{r,o,z,s}=\EE\left[\left(\frac{\partial^3}{\partial x^r\partial x^o\partial x^z}G_{l+1,j}^s(x)\right)^{n}\right],
\]
then we get, due to $3a^2b^2c^2\le a^6+b^6+c^6$ and
\begin{align*}
2\EE\left[XYZU\right]\le & 2\sqrt{\EE\left[X^2\right]}\sqrt[6]{\EE\left[Y^6\right]}\sqrt[6]{\EE\left[Z^6\right]}\sqrt[6]{\EE\left[U^6\right]}\le \EE\left[X^2\right]+\sqrt[3]{\EE\left[Y^6\right]}\sqrt[3]{\EE\left[Z^6\right]}\sqrt[3]{\EE\left[U^6\right]}\\
&\le \EE\left[X^2\right]+\frac{1}{3}\left(\EE\left[Y^6\right]+\EE\left[Z^6\right]+\EE\left[U^6\right]\right),
\end{align*}
as well as~\eqref{smooth:Phi_assump_1} and~\eqref{smooth:Phi_assump_2},
\begin{eqnarray*}
\rho_{l+1,2,3}^{r,o,z,k} & \leq & (1+A_{2}\Delta)\rho_{l,2,3}^{r,o,z,k}
+\sum_{s:\,s\neq k}\left(C_{1,1}\Delta(\rho_{l,2,3}^{r,o,z,k}+\rho_{l,2,3}^{r,o,z,s})+(d-1)B_{2}\Delta\rho_{l,2,3}^{r,o,z,s}\right)\\
&&+\sum_{s,u,v=1}^d C_{1,1}\Delta\left(\rho_{l,2,2}^{r,o,z,v}+\frac{1}{2}\left(\rho_{l,4,1}^{r,s}+\rho_{l,4,1}^{o,u}+\rho_{l,4,1}^{z,u}+\rho_{l,4,2}^{r,o,s}+\rho_{l,4,2}^{r,z,s}+\rho_{l,4,2}^{o,z,u}\right)\right)\\
&&+\sum_{s,u,v,w=1}^d C_{1,1}\Delta\left(\rho_{l,2,2}^{r,o,z,w}+\frac{1}{3}\left(\rho_{l,6,1}^{r,s}+\rho_{l,6,1}^{o,u}+\rho_{l,6,1}^{z,v}\right)\right)\\
&&+3d^2\sum_{s,u=1}^dB_{2}\Delta\left(\rho_{l,4,1}^{r,s}+\rho_{l,4,1}^{o,u}+\rho_{l,4,1}^{z,u}+\rho_{l,4,2}^{r,o,s}+\rho_{l,4,2}^{r,z,s}+\rho_{l,4,2}^{o,z,u}\right)\\
&&+d^3\sum_{s,u,v=1}^dB_{2}\Delta\frac{1}{3}\left(\rho_{l,6,1}^{r,s}+\rho_{l,6,1}^{o,u}+\rho_{l,6,1}^{z,v}\right).
\end{eqnarray*}
Further, denote 
\[
\rho_{l+1,2,3}^{r,o,z}=\sum_{s=1}^d\rho_{l+1,2,3}^{r,o,z,s},
\]
then we get 
\begin{eqnarray*}
\rho_{l+1,2,3}^{r,o,z} & \leq & (1+A_{2}\Delta)\rho_{l,2,2}^{r,o,z}
+2(d-1)C_{1,1}\Delta\rho_{l,2,2}^{r,o,z}+(d-1)^2B_{2}\Delta\rho_{l,2,2}^{r,o,z}\\
&&+d^3C_{1,1}\Delta\left(\rho_{l,2,2}^{r,o,z}+\frac{1}{2}\left(\rho_{l,4,1}^{r}+\rho_{l,4,1}^{o}+\rho_{l,4,1}^z+\rho_{l,4,2}^{r,o}+\rho_{l,4,2}^{r,z}+\rho_{l,4,2}^{o,z}\right)\right)\\
&&+d^4C_{1,1}\Delta\left(\rho_{l,2,2}^{r,o,z}+\frac{1}{3}\left(\rho_{l,6,1}^{r}+\rho_{l,6,1}^{o}+\rho_{l,6,1}^{z}\right)\right)\\
&&+3d^4B_2\Delta\left(\rho_{l,4,1}^{r}+\rho_{l,4,1}^{o}+\rho_{l,4,1}^z+\rho_{l,4,2}^{r,o}+\rho_{l,4,2}^{r,z}+\rho_{l,4,2}^{o,z}\right)\\
&&+d^6B_2\Delta\frac{1}{3}\left(\rho_{l,6,1}^{r}+\rho_{l,6,1}^{o}+\rho_{l,6,1}^{z}\right).
\end{eqnarray*}
This gives us
\[
\rho_{l+1,2,3}^{r,o,z}\leq (1+\kappa_{11}\Delta)\rho_{l,2,2}^{r,o,z}+\kappa_{12}\Delta,
\]
for some constants $\kappa_{11},\kappa_{12}>0$,
leading to
\[
\rho_{l,2,2}^{r,o,z}\leq(1+\kappa_{11}\Delta)^{l-j-1}\rho_{j+1,2,3}^{r,o,z}+\kappa_{13}
,\quad l=j+1,\ldots,J-1,
\]
where $\kappa_{13}>0$ and
\[
\rho_{j+1,2,3}^{r,o,z}=\sum_{s=1}^d\EE\left[\left(\frac{\partial^3}{\partial x^r\partial x^o\partial x^z}\Phi_{\Delta}^s\left(x,\xi_{j+1}\right)\right)^{2}\right].
\]
Thus, we obtain the boundedness of
\begin{eqnarray*}
&&\left|\frac{\partial^3}{\partial x^r\partial x^o\partial x^z}q_j(x)\right|
\\
&\le&\sum_{s=1}^d\EE\left|\frac{\partial }{\partial x^s}f(G_{J,j}(x))\frac{\partial^3}{\partial x^r\partial x^o\partial x^z}G_{J,j}^s(x)\right|\\
&&+\sum_{s,u=1}^d\EE\left|\frac{\partial^2 }{\partial x^s\partial x^u}f(G_{J,j}(x))\left(\frac{\partial^2}{\partial x^r\partial x^o}G_{J,j}^s(x)\frac{\partial}{\partial x^z}G_{J,j}^u(x)+\frac{\partial^2}{\partial x^r\partial x^z}G_{J,j}^s(x)\frac{\partial}{\partial x^o}G_{J,j}^u(x)\right.\right.\\
&&\phantom{+\sum_{s,u=1}^d\EE\left|\frac{\partial^2 }{\partial x^s\partial x^u}f(G_{J,j}(x))\left(\right.\right.}\left.\left.+\frac{\partial}{\partial x^r}G_{J,j}^s(x)\frac{\partial^2}{\partial x^o\partial x^z}G_{J,j}^u(x)\right)\right|\\
&&+\sum_{s,u,v=1}^d\EE\left|\frac{\partial^3 }{\partial x^s\partial x^u\partial x^v}f(G_{J,j}(x))\frac{\partial}{\partial x^r}G_{J,j}^s(x)\frac{\partial}{\partial x^o}G_{J,j}^u(x)\frac{\partial}{\partial x^z}G_{J,j}^v(x)\right]\\
&\le &\sum_{s=1}^d\sqrt{\EE\left[\left(\frac{\partial }{\partial x^s}f(G_{J,j}(x))\right)^2\right]\rho_{J,2,3}^{r,o,z,s}}\\
&&+\sum_{s,u=1}^d\sqrt{\EE\left[\left(\frac{\partial^2 }{\partial x^s\partial x^u}f(G_{J,j}(x))\right)^2\right]}\left(\sqrt[4]{\rho_{J,4,2}^{r,o,s}\rho_{J,4,1}^{z,u}}+\sqrt[4]{\rho_{J,4,2}^{r,z,s}\rho_{J,4,1}^{o,u}}+\sqrt[4]{\rho_{J,4,1}^{r,s}\rho_{J,4,2}^{o,z,u}}\right)\\
&&+\sum_{s,u,v=1}^d\sqrt{\EE\left[\left(\frac{\partial^3 }{\partial x^s\partial x^u\partial x^v}f(G_{J,j}(x))\right)^2\right]}\sqrt[6]{\rho_{J,6,1}^{r,s}\rho_{J,6,1}^{o,u}\rho_{J,6,1}^{z,v}}
\end{eqnarray*}
for all $r,o,z\in\left\{1,\ldots,d\right\}$, since $f$ is assumed to be three times continuously differentiable with bounded partial derivatives up to order~3.

\subsection{Proof of Theorem~\protect\ref{th:first_old}}
(i) Straightforward.

(ii) Let us define $\mu_\Delta(x):=x+\mu(x)\Delta$. Then we obtain via Taylor's theorem (cf.~\eqref{eq:PhiK=1})
\begin{align*}
q_j(\Phi_\Delta(x,y))=q_j(\mu_\Delta(x))+\sqrt{\Delta}\sum_{k=1}^d\sum_{i=1}^m  \sigma^{ki}(x)y^i\int\limits_0^1\frac{\partial q_j}{\partial x^k}(\mu_\Delta(x)+t\sigma(x)\sqrt{\Delta}y)\,dt.
\end{align*}
This gives us (see~\eqref{eq:coef05a})
\begin{align}
\notag
a_{j,r,s}(x)&=\frac{1}{2^m}\sum_{y\in\left\{-1,1\right\}^m}
q_j(\Phi_\Delta(x,y))\prod_{o=1}^r y^{s_o}\\
\label{rel_a_dq}
&=\frac{\sqrt{\Delta}}{2^m}\sum_{y\in\left\{-1,1\right\}^m}\left(\prod_{o=1}^ry^{s_o}\right)\sum_{k=1}^d\sum_{i=1}^m  \sigma^{ki}(x)y^i\int\limits_0^1\frac{\partial q_j}{\partial x^k}(\mu_\Delta(x)+t\sigma(x)\sqrt{\Delta}y)\,dt,
\end{align}
since
$$
\frac{1}{2^m}\sum_{y\in\left\{-1,1\right\}^m}\prod_{o=1}^ry^{s_o}=\EE\left[\prod_{o=1}^r\xi_j^{s_o}\right]=0.
$$
Next we apply Theorem~\ref{smooth:func_assump_general} for the case $K=1$ to get that all the functions $q_j$ are continuously differentiable with bounded partial derivatives. Clearly, the assumptions in this theorem hold, when all the functions $f,\mu^k,\sigma^{ki}$, $k\in\left\{1,\ldots,d\right\}$, $i\in\left\{1,\ldots,m\right\}$ are continuously differentiable with bounded derivatives. Together with the assumption, that all the functions $\sigma^{ki}$ are bounded, we get from~\eqref{rel_a_dq} that $a_{j,r,s}$ is of order $\sqrt{\Delta}$ for all $j,r,s$.

\subsection{Proof of Theorem~\protect\ref{th:first_new}}
Let us consider a higher order Taylor expansion compared to the proof of Theorem~\ref{th:first_old} and recall that $\mu_\Delta(x)=x+\mu(x)\Delta$. We have for any $y\in\left\{-1,1\right\}^m$
\begin{align}
\notag
q_j(\Phi_\Delta(x,y))=&q_j(\mu_\Delta(x))+\sqrt{\Delta}\sum_{k=1}^d\frac{\partial}{\partial x^k}q_j(\mu_\Delta(x))\sum_{i=1}^m\sigma^{ki}(x)y^i\\
\label{out:q_first}
&+\Delta\sum_{k,l=1}^d(2-\delta_{k,l})\int\limits_0^1(1-t)\frac{\partial^2}{\partial x^k\partial x^l}q_j(\mu_\Delta(x)+t\sqrt{\Delta}\sigma(x)y)\,dt\sum_{i=1}^m\sigma^{ki}(x)y^i\sum_{i=1}^m\sigma^{li}(x)y^i,
\end{align}
where $\delta_{\cdot,\cdot}$ is the Kronecker delta. This gives us for $r\ge 2$ (cf.~\eqref{eq:coef05a})
\begin{align}
\notag
a_{j,r,s}(x)&=\frac{1}{2^m}\sum_{y\in\left\{-1,1\right\}^m}
q_j(\Phi_\Delta(x,y))\prod_{o=1}^r y^{s_o}\\
\notag
&=\frac{\Delta}{2^m}\sum_{k,l=1}^d\left[(2-\delta_{k,l})\sum_{y\in\left\{ -1,1\right\} ^{m}}\left(\sum_{i=1}^m\sigma^{ki}(x)y^i\sum_{i=1}^m\sigma^{li}(x)y^i\prod_{o=1}^r y^{s_o}\right.\right.\\
\label{rel_a_dq2}
&\phantom{=\frac{\Delta}{2^m}\sum_{k,l=1}^d\left(\right.}\left.\left.\cdot\int\limits_0^1(1-t)\frac{\partial^2}{\partial x^k\partial x^l}q_j(\mu_\Delta(x)+t\sqrt{\Delta}\sigma(x)y)\,dt\right)\right],
\end{align}
due to (cf.~\eqref{out:q_first})
\begin{align}
\label{out:first_r2}
\frac{1}{2^m}\sum_{y\in\left\{-1,1\right\}^m}y^i\prod_{o=1}^r y^{s_o}=\EE\left[\xi_j^i\prod_{o=1}^r \xi_j^{s_o}\right]=0
\end{align}
for all $i\in\left\{1,\ldots,m\right\}$. (Note that~\eqref{out:first_r2} does not hold for $r=1$.) Applying Theorem~\ref{smooth:func_assump_general} (case $K=2$), we get that $q_j$ is twice continuously differentiable with bounded partial derivatives up to order 2, provided that all the functions $f,\mu^k,\sigma^{k,i}$ are twice continuously differentiable with bounded partial derivatives up to order 2. Together with the assumption, that all the functions $\sigma^{kl}$ are bounded, we get from~\eqref{rel_a_dq2} that $a_{j,r,s}$ is of order $\Delta$ for all $j,r,s$ with $r>1$.

\subsection{Proof of Theorem~\protect\ref{th:first_var}}
Here we apply Theorem~\ref{th:first_new}, which gives us (cf.~\eqref{eq:2909a2})
\begin{align*}
\Var\left[f(X_{\Delta,T})-M_{\Delta,T}^{(1),trunc}\right]&=\Var\left[M_{\Delta,T}^{(1)}-M_{\Delta,T}^{(1),trunc}\right]\\
&=\Var\left[\sum_{j=1}^J
\sum_{r=2}^m
\sum_{1\le s_1<\ldots<s_r\le m}
a_{j,r,s}(X_{\Delta,(j-1)\Delta})
\prod_{i=1}^r \xi_j^{s_i}\right] \\
&\lesssim\sum_{j=1}^J
\sum_{r=2}^m
\sum_{1\le s_1<\ldots<s_r\le m}\EE\left[a_{j,r,s}^2(X_{\Delta,(j-1)\Delta})
\right]\lesssim \Delta,
\end{align*}
since $\EE\left[a_{j,r,s}^2(X_{\Delta,(j-1)\Delta})
\right]\lesssim\Delta^2$ for all $j,r,s$ with $r>1$.

\subsection{Proof of Theorem~\protect\ref{th:second_old}}
The proof works similarly to the one of Theorem~\ref{th:first_old}. More precisely, here we define (cf.~\eqref{eq:2002a6})
$$
\mu_\Delta(x):=x+\mu(x)\Delta+\frac12\cL^0\mu(x)\,\Delta^2.
$$
Then we derive the zero-order Taylor expansion for $q_j(\Phi_\Delta(x,y,z))$ around $\mu_\Delta(x)$, use that 
\[
\EE\left[\prod_{r\in U_1} H_{o_r}(\xi_j^r)
\prod_{(k,l)\in U_2} V_j^{kl}\right]=0
\]
 and observe that all components $\tilde\Phi_\Delta^k(x,y,z):=\Phi_\Delta^k(x,y,z)-\mu_\Delta^k(x)$, $k\in\left\{1,\ldots,d\right\}$ (as an analogue of $\sqrt{\Delta}\sum_{i=1}^m\sigma^{ki}(x)y^i$ in case of the weak Euler scheme), are of order $\sqrt{\Delta}$ under less strict assumptions than required in the present theorem. Finally we apply Theorem~\ref{smooth:func_assump_general} (case $K=1$) which gives us that $q_j$ is continuously differentiable with bounded partial derivatives under the assumptions, that all functions $\mu^k$ and $\sigma^{ki}$ are bounded and all the functions $f,\mu^k,\sigma^{ki}$ are three-times continuously differentiable with bounded partial derivatives up to order 3. Consequently, all the functions $a_{j,o,U_1,U_2}$ are of order $\sqrt{\Delta}$.

\subsection{Proof of Theorem~\protect\ref{th:second_new}}
(i) The proof works similarly to the one of Theorem~\ref{th:first_new}, that is, we consider a Taylor expansion for $q_j(\Phi_\Delta(x,y,z))$ of order 1, around the same point $\mu_\Delta(x)$ as in the proof of Theorem~\ref{th:second_old}. Then we use 
$$\EE\left[\tilde\Phi_\Delta^k(x,\xi_j,V_j)\prod_{r\in U_1} H_{o_r}(\xi_j^r)
\prod_{(k,l)\in U_2} V_j^{kl}\right]=0,\quad k\in\left\{1,\ldots,d\right\},$$
whenever $\left|U_2\right|+\left|\mathcal{K}_2\right|+\frac{\left|\mathcal{K}_1\right|}{2}\ge 1$ (where again $\tilde\Phi_\Delta^k(x,y,z)=\Phi_\Delta^k(x,y,z)-\mu_\Delta^k(x)$). Then
we apply Theorem~\ref{smooth:func_assump_general} (case $K=2$) which gives us that $q_j$ is twice continuously differentiable with bounded partial derivatives up to order 2 under the assumptions, that all functions $\mu^k$ and $\sigma^{ki}$ are bounded and all the functions $f,\mu^k,\sigma^{ki}$ are four-times continuously differentiable with bounded partial derivatives up to order 4. Finally, we get that all the functions $a_{j,o,U_1,U_2}$ are of order $\Delta$, since the product of all functions $\tilde\Phi_\Delta^k(x,y,z)\tilde\Phi_\Delta^l(x,y,z)$, $k,l\in\left\{1,\ldots,d\right\}$, is of order $\Delta$ under the above assumptions.

(ii) Here we consider the Taylor expansion of order 2, that is
\begin{align*}
q_j(\Phi_\Delta(x,y,z))=&q_j(\mu_\Delta(x))+\sum_{k=1}^d\frac{\partial}{\partial x^k}q_j(\mu_\Delta(x))\tilde\Phi_\Delta^k(x,y,z)\\
&+\sum_{k,l=1}^d\frac{1}{2}(2-\delta_{k,l})\frac{\partial^2}{\partial x^k\partial x^l}q_j(\mu_\Delta(x))\tilde\Phi_\Delta^k(x,y,z)\tilde\Phi_\Delta^l(x,y,z)\\
&+\sum_{k,l,n=1}^d\left[\left(3-\frac{3}{2}\left(\delta_{k,l}+\delta_{k,n}+\delta_{l,n}\right)+2\delta_{k,l}\delta_{k,n}\delta_{l,n}\right)\tilde\Phi_\Delta^k(x,y,z)\Phi_\Delta^l(x,y,z)\Phi_\Delta^n(x,y,z)\right.\\
&\phantom{+\sum_{k,l,n=1}^d\left[\right.}\left.
\cdot\int\limits_0^1(1-t)^2\frac{\partial^3}{\partial x^k\partial x^l\partial x^n}q_j(\mu_\Delta(x)+t\tilde\Phi_\Delta(x,y,z))\,dt\right].
\end{align*}
Next we use 
$$\EE\left[\tilde\Phi_\Delta^k(x,\xi_j,V_j)\tilde\Phi_\Delta^l(x,\xi_j,V_j)\prod_{r\in U_1} H_{o_r}(\xi_j^r)
\prod_{(k,l)\in U_2} V_j^{kl}\right]=0,\quad k,l\in\left\{1,\ldots,d\right\},$$
whenever $\left|U_2\right|+\left|\mathcal{K}_2\right|+\frac{\left|\mathcal{K}_1\right|}{2}>1$, and thus we obtain (cf.~\eqref{eq:0403a11})
\begin{align*}
&a_{j,o,U_1,U_2}(x)\\
=&\frac{1}{2^{\frac{m(m-1)}{2}}}\,\frac{1}{6^{m}}\sum_{y\in\{-\sqrt{3},0,\sqrt{3}\}^{m}}\sum_{z\in\{-1,1\}^{\frac{m(m-1)}{2}}} 4^{\sum_{i=1}^{m}I(y^{i}=0)}
\prod_{r\in U_1}H_{o_r}(y^r)
\prod_{(k,l)\in U_2}z^{kl}\,
q_j(\Phi_\Delta(x,y,z))\\
=&\frac{1}{2^{\frac{m(m-1)}{2}}}\,\frac{1}{6^{m}}\sum_{y\in\{-\sqrt{3},0,\sqrt{3}\}^{m}}\sum_{z\in\{-1,1\}^{\frac{m(m-1)}{2}}} 4^{\sum_{i=1}^{m}I(y^{i}=0)}
\prod_{r\in U_1}H_{o_r}(y^r)
\prod_{(k,l)\in U_2}z^{kl}\\
&\cdot\sum_{k,l,n=1}^d\left[\left(3-\frac{3}{2}\left(\delta_{k,l}+\delta_{k,n}+\delta_{l,n}\right)+2\delta_{k,l}\delta_{k,n}\delta_{l,n}\right)\tilde\Phi_\Delta^k(x,y,z)\Phi_\Delta^l(x,y,z)\Phi_\Delta^n(x,y,z)\right.\\
&\phantom{+\sum_{k,l,n=1}^d\left[\right.}\left.
\cdot\int\limits_0^1(1-t)^2\frac{\partial^3}{\partial x^k\partial x^l\partial x^n}q_j(\mu_\Delta(x)+t\tilde\Phi_\Delta(x,y,z))\,dt\right]
\end{align*}
Then we apply Theorem~\ref{smooth:func_assump_general} (case $K=3$) which gives us that $q_j$ is three-times continuously differentiable with bounded partial derivatives up to order 3 under the assumptions, that all functions $\mu^k$ and $\sigma^{ki}$ are bounded and all the functions $f,\mu^k,\sigma^{ki}$ are five-times continuously differentiable with bounded partial derivatives up to order 5. Finally, we get that all the functions $a_{j,o,U_1,U_2}$ are of order $\Delta^{3/2}$, since the product of all functions $\tilde\Phi_\Delta^k(x,y,z)\tilde\Phi_\Delta^l(x,y,z)\tilde\Phi_\Delta^n(x,y,z)$, $k,l,n\in\left\{1,\ldots,d\right\}$, is of order $\Delta^{3/2}$ under the above assumptions.

\subsection{Proof of Theorem~\protect\ref{th:second_var}}
The proof is similar to the one of Theorem~\ref{th:first_var}.

\subsection{Proof of Lemma~\protect\ref{th:2104a1}}
We refer to Theorem~11.3 in~\cite{gyorfi2002distribution}. When applying it, we obtain actually
\begin{align}
\label{regr:gjorfi_direct}
\E\|\tilde a_{j,o,U_1,U_2}-a_{j,o,U_1,U_2}\|^2_{L^2(\PP_{\Delta,j-1})}
\le\tilde{c}\max\left\{\Sigma,\tilde A^2\Delta_{U_1,U_2}^2\right\}(\log N+1)\frac{Q}{N}
+\frac{8\,C_\kappa}{Q^\kappa}.
\end{align}
However, the maximum in~\eqref{regr:gjorfi_direct} is in fact a sum of two terms $\Sigma$ and $\tilde A^2\Delta_{U_1,U_2}(\log N_r+1)$ so that the logarithm is only included in one term (see the proof of Theorem~11.3 in~\cite{gyorfi2002distribution}).

\subsection{Proof of Theorem~\protect\ref{th:2106a1}}
Using the martingale transform structure
in~\eqref{eq:28042016a1} and~\eqref{cv:second_trunc}
(recall footnote~\ref{ft:19102016a1}
on page~\pageref{ft:19102016a1}) together with the orthonormality of the system
$\prod_{r\in U_1} H_{o_r}(\xi_j^r)
\prod_{(k,l)\in U_2} V_j^{kl}$, we get by~\eqref{second:cv_var} and~\eqref{eq:2104a2}
\begin{align*}
\Var\left[f(X_{\Delta,T})-
\widetilde M^{(2),trunc}_{\Delta,T}\right]&=\Var\left[f(X_{\Delta,T})-
M^{(2),trunc}_{\Delta,T}\right]+\Var\left[M^{(2),trunc}_{\Delta,T}-
\widetilde M^{(2),trunc}_{\Delta,T}\right]\\
&\lesssim\Delta^2+\sum_{j=1}^J
\sum_{\substack{(U_1,U_2)\in\cA\\ \left|U_2\right|+\left|\mathcal{K}_2\right|+\frac{1}{2}\left|\mathcal{K}_1\right|\le 1}}
\sum_{o\in\{1,2\}^{U_1}}
\EE\|\tilde a_{j,o,U_1,U_2}-a_{j,o,U_1,U_2}\|^2_{L^2(\PP_{\Delta,j-1})}
\\
&\le
\Delta^2+Jm(m+1)\left(\tilde{c}(\Sigma+\tilde A^2\Delta(\log N+1))\frac{Q}{N}
+\frac{8\,C_\kappa}{Q^\kappa}\right),
\end{align*}
since $\Delta_{U_1,U_2}^2\le \Delta$.

\subsection{Proof of Theorem~\protect\ref{trcv_compl}}
The proof is similar to the complexity analysis
performed in~\cite{belomestny2016konakov}.

\bibliographystyle{abbrv}
\bibliography{ml_vr_bibliography}

\begin{thebibliography}{10}

\bibitem{belomestny2016variance}
D.~Belomestny, S.~H\"afner, T.~Nagapetyan, and M.~Urusov.
\newblock Variance reduction for discretised diffusions via regression.
\newblock {\em Preprint, arXiv:1510.03141v3}, 2016.

\bibitem{belomestny2016dual}
D.~Belomestny, S.~H\"afner, and M.~Urusov.
\newblock Regression-based complexity reduction of the dual nested {M}onte
  {C}arlo methods.
\newblock {\em Preprint, arXiv:1611.06344}, 2016.

\bibitem{belomestny2016konakov}
D.~Belomestny, S.~H\"afner, and M.~Urusov.
\newblock Regression-based variance reduction approach for strong approximation
  schemes.
\newblock {\em Preprint, arXiv:1612.03407v2}, 2017.

\bibitem{belomestny2016stratified}
D.~Belomestny, S.~H\"afner, and M.~Urusov.
\newblock Stratified regression-based variance reduction approach for weak
  approximation schemes.
\newblock {\em Preprint,  arXiv:1612.05255v2}, 2017.

\bibitem{giles2008multilevel}
M.~B. Giles.
\newblock Multilevel {M}onte {C}arlo path simulation.
\newblock {\em Operations Research}, 56(3):607--617, 2008.

\bibitem{gyorfi2002distribution}
L.~Gy{\"o}rfi, M.~Kohler, A.~Krzy{\.z}ak, and H.~Walk.
\newblock {\em A distribution-free theory of nonparametric regression}.
\newblock Springer Series in Statistics. Springer-Verlag, New York, 2002.

\bibitem{KP}
P.~Kloeden and E.~Platen.
\newblock {\em Numerical solution of stochastic differential equations},
  volume~23.
\newblock Springer, 1992.

\bibitem{milstein2009practical}
G.~N. Milstein and M.~V. Tretyakov.
\newblock Practical variance reduction via regression for simulating
  diffusions.
\newblock {\em SIAM Journal on Numerical Analysis}, 47(2):887--910, 2009.

\bibitem{muller2015complexity}
T.~M{\"u}ller-Gronbach, K.~Ritter, and L.~Yaroslavtseva.
\newblock On the complexity of computing quadrature formulas for marginal
  distributions of {SDE}s.
\newblock {\em Journal of Complexity}, 31(1):110--145, 2015.

\bibitem{muller1324deterministic}
T.~M{\"u}ller-Gronbach and L.~Yaroslavtseva.
\newblock Deterministic quadrature formulas for {SDE}s based on simplified weak
  {I}t{\^o}-{T}aylor steps.
\newblock {\em Foundations of Computational Mathematics}, 16(5):1325--1366,
  2016.

\bibitem{newton1994variance}
N.~J. Newton.
\newblock Variance reduction for simulated diffusions.
\newblock {\em SIAM journal on applied mathematics}, 54(6):1780--1805, 1994.

\end{thebibliography}
\end{document}